\documentclass[11pt]{amsart}

\usepackage{amssymb}
\usepackage{latexsym}
\usepackage{amsmath}
\usepackage{amsthm}
\usepackage{amsfonts}
\usepackage{color}
\usepackage{pdfsync}
\usepackage[usenames,dvipsnames]{pstricks}
\usepackage{epsfig}
\usepackage{pst-grad} 
\usepackage{pst-plot} 

\newcommand{\R}{\mathbb{R}}
\newcommand{\N}{\mathbb{N}}
\newcommand{\Z}{\mathbb{Z}}
\newcommand{\Q}{\mathcal{Q}}

\newcommand{\I}{\mathcal{I}}

\newcommand{\T}{\mathbb{T}}

\theoremstyle{plain}
\newtheorem{defi}{Definition}[section]
\newtheorem{prop}[defi]{Proposition}
\newtheorem{teo}[defi]{Theorem}
\newtheorem{cor}[defi]{Corollary}
\newtheorem{lema}[defi]{Lemma}

\newtheorem{remark}[defi]{Remark}

\theoremstyle{definition}

\theoremstyle{remark}

\numberwithin{equation}{section}

\begin{document}

\title[]{Regularity Results and Large Time Behavior for Integro-Differential Equations with Coercive Hamiltonians.}

\author[]{Guy Barles}
\address{
Guy Barles:
Laboratoire de Math\'ematiques et Physique Th\'eorique (UMR CNRS 7350), F\'ed\'eration Denis Poisson (FR CNRS 2964),
Universit\'e Fran\c{c}ois Rabelais Tours, Parc de Grandmont, 37200 Tours, FRANCE. {\tt Guy.Barles@lmpt.univ-tours.fr}}

\author[]{Shigeaki Koike}
\address{
Shigeaki Koike: Mathematical Institute, Tohoku University, Aoba, Sendai 980-8578, JAPAN. {\tt koike@math.tohoku.ac.jp}}

\author[]{Olivier Ley}
\address{ 
Olivier Ley: IRMAR, INSA de Rennes, 35708 Rennes, FRANCE. {\tt  olivier.ley@insa-rennes.fr}}

\author[]{Erwin Topp}
\address{
Erwin Topp:
Departamento de Ingenier\'\i a Matem\'atica (UMI 2807 CNRS), Universidad de Chile, Casilla 170, Correo 3, Santiago, CHILE.
and Laboratoire de Math\'ematiques et Physique Th\'eorique (CNRS UMR 7350), Universit\'e Fran\c{c}ois Rabelais,
Parc de Grandmont, 37200, Tours, FRANCE. 
{\tt etopp@dim.uchile.cl}} 

\date{\today}

\begin{abstract} 
In this paper we obtain regularity results for elliptic integro-differential equations driven by the stronger effect of coercive gradient terms. 
This feature allows us to construct suitable strict supersolutions from which we conclude H\"older estimates for bounded subsolutions. In many 
interesting situations, this gives way to a priori estimates for subsolutions. 
We apply this regularity results to obtain the ergodic asymptotic behavior of the associated 
evolution problem in the case of superlinear equations. One of the surprising features in our proof 
is that it avoids the key ingredient which are usually necessary to use the Strong Maximum Principle: 
linearization based on the Lipschitz regularity of the solution of the ergodic problem. The proof entirely relies 
on the H\"older regularity.
\end{abstract}



\maketitle


\section{Introduction.}

In~\cite{Capuzzo-Dolcetta-Leoni-Porretta}, Capuzzo-Dolcetta, Leoni and Porretta prove a surprising regularity result for {\em subsolutions} of 
superquadratic second-order elliptic equations which can be described in the following way. We consider the model equation
\begin{equation}\label{superquadratic}
\lambda v - Tr(A(x) D^2v(x)) + b(x) |Dv(x)|^m = f(x) \quad \hbox{in   }\Omega, 
\end{equation}
where $ \Omega$ is an open subset of $\R^N$, $A, b,f$ are continuous functions in $\Omega$, $A$ taking values in the set of 
nonnegative matrices and $b, f$ are real valued, with $b(x)\geq b_0 >0$ in $\Omega$, $m > 2$ and $\lambda \geq 0$. 
The function $v:\Omega \to \R$ is a real-valued solution and $Dv,D^2v$ denote 
its gradient and Hessian matrix. 
In~\cite{Capuzzo-Dolcetta-Leoni-Porretta}, the authors prove that, if $u : \Omega \to \R$ is a bounded viscosity subsolution of (\ref{superquadratic}) 
then $u$ is locally H\"older continuous with exponent $\alpha := (m-2)(m-1)^{-1}$ and the local H\"older seminorm depends only on the datum
($L^\infty$ bounds on $A, f$ and $b_0$) but not on any $L^\infty$ bound nor oscillation of $u$. Actually this result provides, 
in many interesting situations, an estimate on the $L^\infty$ norm of $u$.

The starting point of the present work was to investigate how such a result could be extended to the case of nonlocal elliptic equations like 
\begin{equation}\label{eq}\tag{P}
\lambda u(x) -I_x(u,x) + H(x, Du(x)) = 0 \quad \hbox{in   }\Omega,
\end{equation}
where $\lambda \geq 0$ and $H : \Omega\times \R^N \to \R$  is a continuous 
nonlinearity having the same properties as $b(x)|p|^m -f(x)$ above. 
The term $I_x$ is a nonlocal operator playing the role of the diffusion, defined as follows:
for $x, y \in \R^N$ and $\phi : \R^N \to \R$ a bounded continuous function which is 
$C^2$ in a neighborhood of $y$, we write
\begin{equation}\label{operator}
I_x(\phi, y) = \int_{\R^N} [\phi(y + z) - \phi(y) - \mathbf{1}_B \langle D\phi(y),  z \rangle] \nu_x(dz),
\end{equation}
where $B$ denotes the unit ball and $\{ \nu_x \}_{x \in \R^N}$ is a family of \textsl{L\'evy measures}, see (M1)-(M2) below for precise assumptions.
An important example of such nonlocal operator is the case when $\nu_x = \nu$ for all $x \in \R^N$, with
\begin{equation*}
\nu(dz) = C_{N, \sigma} |z|^{-(N + \sigma)}dz, 
\end{equation*}
where $\sigma \in (0,2)$ and $C_{N,\sigma}$ is a normalizing constant. In that case, for all $x \in \R^N$, $-I_x = (-\Delta)^{\sigma/2}$ is
the fractional Laplacian of order $\sigma$ (see~\cite{Hitch}).
By the form of $I_x$ in~\eqref{operator}, we point out that subsolutions of~\eqref{eq} must be defined on $\R^N$ or at least 
in a large enough domain (depending on $\nu_x$) in order that the nonlocal operator is well-defined. 

In \cite{Capuzzo-Dolcetta-Leoni-Porretta} and even more in the simplified version given in \cite{Barles1}, 
the authors take advantage of the superquadratic gradient term to construct {\em locally} a strict supersolution to~\eqref{superquadratic} using 
power-like functions. The power profile of such supersolutions gives the (local) H\"older regularity for 
bounded \textsl{subsolutions} of the equation.
This proof is based on the leading  effect of the gradient term more than on the ellipticity, resembling the behavior of first-order coercive 
equations (see~\cite{Barles-book}). The H\"older exponent $(m - 2)(m - 1)^{-1}$ just comes from a simple balance of powers in~\eqref{superquadratic} 
and this H\"older regularity can be extended up to the boundary of the domain if it is regular enough (see also~\cite{Barles1}). 

All these arguments seem extendable to the nonlocal framework, and in particular, if we think the nonlocal term as an operator of order $\sigma \in (0,2)$. 
But here is a key difference which is going to play a double role : first, depending on the support of the measure $\nu_x(dz)$, the operator may use values of $u$ outside $\Omega$. This arises, typically, when equation~\eqref{eq} is complemented by an exterior Dirichlet condition (see~\cite{Barles-Chasseigne-Imbert}). Of course, and this is very natural in the case of exterior Dirichlet condition, these outside values cannot be controlled by the equation. Hence, in that case, it is clearly impossible to have results which are independent of the $L^\infty$ norm or oscillation of $u$.

On the contrary, this analysis shows that, in principle, this could be possible in the case when the support of the measure is such that the integral 
of $I_x(u,y)$ only takes into account points such that $y+z\in \Omega$, typically when
\begin{equation}\label{Cens-FL}
I_x(\phi, y) = C_{N, \sigma} \int_{y+z \in \Omega} [\phi(y + z) - \phi(y) - \mathbf{1}_B \langle D\phi(y),  z \rangle] |z|^{-(N + \sigma)}dz\; .
\end{equation}

These type of operators are related to ``censored processes'' in the probabilistic literature : in this context, it means that the jumps 
processes cannot jump from $\Omega$ to $\Omega^c$. We refer to
e.g. \cite{BBC,FOT,GQY, GQY2,KP,NJ:Book} for more details on such processes. In \cite{BBC,GQY2}, the censored fractional Laplacian appears in 
connection with Dirichlet forms; they also appear in the analysis literature as regional Laplacians (\cite{Ishii-Nakamura}) and very 
naturally in the study of Neumann boundary conditions (\cite{bcgj}). We therefore call \textsl{censored operators (with respect to $\Omega$)}
the operators which satisfy
\begin{equation*}
x + \mathrm{supp} \{ \nu_x \} \subset \Omega, \quad \mbox{for all} \ x \in \Omega.
\end{equation*}

Actually, we remark that we can always reduce to the case of a censored operator by incorporating the integral over the complement
of $\Omega$ into the right-hand side $f$ (see Lemma~\ref{lemacensored} below and/or~\cite{Topp}). 
This ``censoring'' procedure modifies 
the right-hand side into a function which blows up at the boundary of $\Omega$ with a rate which is controlled in 
terms of the singularity of the measure (the $\sigma$ in the fractional Laplacian case) and the oscillation of $u$. Thus, as it can be 
seen in~\cite{Capuzzo-Dolcetta-Leoni-Porretta}, the presence of these unbounded ingredients in the equation restricts
the expected values of the H\"older exponent if we wish a result which holds up to the boundary. Moreover, the same effect arises even for
nonlocal operators
which are originally censored, since the proof of the H\"older regularity consists in localizing, 
typically in some ball included in $\Omega$ and, at this step also, the values of $u$ outside the ball creates essentially the same difficulty 
as the one described above : if we want to write the nonlocal equation as a censored equation in the ball, then this mechanically changes 
the ``natural'' H\"older exponent because of the right-hand side which blows up at the boundary of the ball. 

All these difficulties explain all the different formulations we give for some results but also the nature of the H\"older exponent we obtain.
To be more specific, we consider the basic model equation 
\begin{equation}\label{eqmodel}
\lambda u(x) + a(x) (-\Delta)^{\sigma/2} u(x) + b(x) |Du(x)|^m = f(x) \quad \hbox{in  }\Omega,
\end{equation}
where $\lambda, b, f$ are as in~\eqref{superquadratic} and $a$ is a continuous real-valued function with $a \geq 0$
in $\Omega$. The role of the superquadraticity in~\eqref{superquadratic} is played by 
a \textsl{superfractional growth condition on the gradient}, which is encoded by $m$ in~\eqref{eqmodel} 
through the assumption
\begin{equation}\label{superfractional}
m > \sigma,
\end{equation}
and the strict positivity requirement on $b$.
The difficulties we mention above on the nonlocality have a price and this price is a ``less natural'' H\"older 
exponent $(m - \sigma)/m$ for subsolutions to~\eqref{eqmodel}. 
Nevertheless, we can get interior H\"older regularity results with ``more natural'' exponents 
$(m - \sigma)/(m - 1)$ if $\sigma > 1$, Lipschitz continuity 
if $\sigma < 1$, and any exponent in $(0,1)$ for $\sigma = 1$, 
since localization arguments are unnecessary in this situation. Finally, we point out that in the case of censored operators 
(here if $(-\Delta)^{\sigma/2}$ is replaced by the operator given by \eqref{Cens-FL}), 
we recover a complete control on the oscillation of $u$ on $\overline{\Omega}$ as a consequence of the form of the estimates 
(see Corollary~\ref{corosc} below).

It is worth pointing out that our results share (with some limitations we described above) the same interesting consequence as the ones of~\cite{Capuzzo-Dolcetta-Leoni-Porretta}, namely a control on the oscillation of 
(sub)solutions to~\eqref{superquadratic} inside $\Omega$ (i.e. at least locally) which is stable as $\lambda \to 0^+.$ This feature has important applications on the study of 
large time behavior for associated parabolic problems and homogenization because of the importance of the ergodic problem.

We are able to provide {\em global} oscillation bounds satisfying this stability property for some class of problems~\eqref{eq} as, for example, 
equations associated to censored operators and obviously for equations set in the whole space $\R^N$. 
This contrasts with the results obtained by Cardaliaguet and Rainer~\cite{Cardaliaguet-Rainer} (see also~\cite{Cardaliaguet-Cannarsa}), 
where the authors obtain very interesting regularity results for (parabolic) superquadratic integro-differential equations using 
a probabilistic approach, but where their H\"older estimates depend on the $L^\infty$ norm of the solution. 

In the second part of this paper, we present an application of our regularity results to the 
study of the large time behavior for Cauchy problems 
\begin{equation}\label{pareq}\tag{CP}
\partial_t u(x, t) - I_x(u(\cdot, t), x) + H(x, Du(x,t)) = 0 \quad \hbox{in   }Q,
\end{equation}
where $Q = \R^N \times (0,+\infty)$.
The asymptotic behavior of the nonlocal evolution problem is also motivated
by its second-order parallel, as the model equation
\begin{equation}\label{evolutionlocal}
\partial_t u(x,t) - \mathrm{Tr}(A(x) D^2u(x,t)) + b(x) |Du(x,t)|^m = f(x) \quad \hbox{in   }Q.
\end{equation}

In the superquadratic case $m> 2$, this evolution equation is also influenced by the stronger 
effect of the first-order term. This can be seen in the paper of Barles and Souganidis~\cite{Barles-Souganidis}, where the authors 
study general equations including~\eqref{superquadratic} and~\eqref{evolutionlocal},
obtain Lipschitz bounds for the solutions
 and prove that, in the periodic setting, the solution 
approaches to the solution of the so-called \textsl{ergodic problem} as $t \to +\infty$. 
This ergodic problem is solved by passing to the limit as $\lambda \to 0^+$ in equation~\eqref{superquadratic},
which is possible by the compactness given by the Lipschitz bounds which are
independent of $\lambda.$
A second key ingredient in the analysis
of the ergodic problem and the large time behavior of~\eqref{evolutionlocal} is the Strong Maximum Principle (\cite{Bardi-DaLio}).

Similar methods and results to~\cite{Barles-Souganidis} are obtained in~\cite{Tchamba} in the context of Cauchy-Dirichlet 
second-order evolution problems in bounded domains. 
In the nonlocal context, analogous ergodic large time behavior for evolution problems are available. 
For instance, in~\cite{Barles-Chasseigne-Ciomaga-Imbert} the authors follow the arguments of~\cite{Barles-Souganidis}, using the Lipschitz regularity results given in~\cite{Barles-Chasseigne-Imbert-Ciomaga-lip}, which allows to ``linearize'' the equation in order to apply the Strong Maximum Principle of~\cite{Ciomaga}. 

In this paper we also follow the lines of~\cite{Barles-Souganidis} to prove the ergodic asymptotic behavior. 
However, contrarily to \cite{Barles-Souganidis} or~\cite{Barles-Chasseigne-Ciomaga-Imbert}, we do not use the Strong Maximum Principle 
in the same way : we do not perform any ``linearization'' of the equation 
(which would have required Lipschitz bounds) and therefore we are able 
to provide results which just use the H\"older regularity of the solutions. This proof requires slightly stronger assumptions on the nonlocal 
operator since we have to be able to use the Strong Maximum Principle \`a 
la Coville~\cite{Coville, Coville2} and to do so, we need the support of the measure defining the nonlocal operator to satisfy an ``iterative covering property''. Though a restriction, this property allows us to study the large time behavior 
for equations associated to very degenerate $x$-dependent nonlocal operators
and $x$-dependent Hamiltonians with a higher degree of coercivity. 

Of course, comparison principles are 
of main importance in this method and for this reason we should focus on a particular class
of $x$-dependent nonlocal operators in \textsl{L\'evy-Ito} form (see~\eqref{operatorLI}). We refer to~\cite{Barles-Imbert} for comparison results 
associated to these operators.

The paper is organized as follows: Section~\ref{regularitysection} is entirely devoted to the regularity results for the stationary problem.
In section~\ref{comparisonsection} we provide the comparison principle and well-possedness of the evolution problem. 
Finally, the large time behavior for this problem is presented in section~\ref{sectionLTB}, where 
the mentioned version
of the strong maximum principle is established.


\medskip
\noindent
{\bf Basic Notation.}
For $x \in \R^N$ and $r > 0$, we denote $B_r(x)$ as the open ball centered at $x$ with radius $r$. We just write $B_r$ for $B_r(0)$
and $B$ for $B_1(0)$.

Let $\Omega \subset \R^N$. We denote as $d_\Omega$ the signed distance function to $\partial \Omega$ which is nonnegative in $\bar{\Omega}$. 
For $\delta > 0$, we also denote $\Omega^\delta = \{x \in \Omega : d_\Omega(x) > \delta\}$.
For any $u:\Omega\to \R,$ the oscillation of $u$ over $\Omega$ is defined by
$${\rm osc}_{\Omega} u= \mathop{\rm sup}_{\Omega} u -  \mathop{\rm inf}_{\Omega} u.$$

For $x, \xi, p \in \R^N$, $A \subset \R^N$ and $\phi$ a bounded function, we define
\begin{equation}\label{defIdelta}
I_\xi[A](\phi, x, p) = \int_{\R^N \cap A} [\phi(x + z) - \phi(x) - \mathbf{1}_{B} \langle p, z \rangle] \nu_\xi(dz).
\end{equation}

We write in a simpler way $I_\xi[A](\phi, x) = I_\xi[A](\phi, x, D\phi(x))$ when $\phi \in L^\infty(\R^N) \cap C^2(B_\delta)$ for some $\delta > 0$,
$I_\xi(\phi, x, p) = I_\xi[\R^N](\phi, x, p)$ when $A = \R^N$ and $I = I_\xi$ if $\nu_\xi = \nu$ does not depend on $\xi$. Note that with these 
notations, $I_x(\phi, x) = I_x[\R^N](\phi, x, D\phi(x))$ for $\phi$ bounded and smooth at $x$ (see~\eqref{operator}).

This paper is based on the viscosity theory to get the results. We refer to~\cite{usersguide, Barles-book, Beginners-Koike}
for the definition and main results of the classical theory, and 
to~\cite{Barles-Imbert, Barles-Chasseigne-Imbert, Alvarez-Tourin, Sayah1, Sayah2} for the nonlocal setting.
Following the definition introduced in the mentioned references, we always assume a viscosity subsolution is upper semicontinuous and 
a viscosity supersolution is lower semicontinuous in the set where the equation takes place.


\section{Regularity.}
\label{regularitysection}

\subsection{Assumptions and Main Regularity Results.}
Let $\sigma \in (0,2)$ fixed. Recalling $I_x$ defined
in~\eqref{operator}, we assume the following conditions over the family $\{ \nu_x \}_x$

\medskip
\noindent
{\bf (M1)} \textsl{ For all $R > 0$ and $\alpha \in [0,2]$, there exists a constant $C_R > 0$ such that, for all $\delta > 0$ we have
\begin{equation*}
\sup \limits_{x \in \bar{B}_R} \int_{B_\delta^c} \min \{1, |z|^\alpha \} \nu_x(dz) 
\leq C_R h_{\alpha, \sigma}(\delta),
\end{equation*}
where $h_{\alpha, \sigma}(\delta)$ is defined for $\delta > 0$ as
\begin{equation}\label{defh}
h_{\alpha, \sigma}(\delta) = \left \{ \begin{array}{cl} \delta^{\alpha - \sigma} \quad & \mbox{if} \ \alpha < \sigma \\
|\ln(\delta)| + 1 \quad & \mbox{if} \ \alpha = \sigma \\
1 \quad & \mbox{if} \ \alpha > \sigma, \end{array} \right . 
\end{equation}
and where we use the convention $|z|^\alpha = 1, z \in \R^N$ when $\alpha = 0$.
}

\medskip
\noindent
{\bf (M2)} \textsl{For all $R > 0$ and $\alpha \in (\sigma, 2]$ there exists a constant $C_R > 0$ such that, 
for all $\delta \in (0,1)$ we have
\begin{equation*}
\sup \limits_{x \in \bar{B}_R} \int_{B_\delta} |z|^\alpha \nu_x(dz) \leq C_R \delta^{\alpha - \sigma}. 
\end{equation*}
}

Assumptions (M1) and (M2) say the nonlocal operator $I_x$ is at most of order $\sigma$, locally in $x \in \R^N$. 
Concerning this last fact, we remark that in the case $\nu_x$ is symmetric and $\sigma \in (0,1)$, $I_x$ defined in~\eqref{operator}
can be written as
\begin{equation}\label{operatorsigmasmall}
I_x(\phi, y) = \int_{\R^N} [\phi(y + z) - \phi(y)]\nu_x(dz),
\end{equation}
for all $y \in \R^N$ and $\phi$ bounded and $C^1$ in a neighborhood of $y$. Since our interest is to keep $I_x$ as a nonlocal operator of 
order $\sigma$, we adopt this formula as a definition for $I_x$ in the case $\sigma \in (0,1)$, even if $\nu_x$ is not symmetric.

In order to expand the application of our results, we consider an open set
$\Omega \subseteq \R^N$ not necessarily bounded, and $H$ satisfying the growth condition
\begin{equation}\label{H}
H(x, p) \geq b_0 |p|^m  - A(d_\Omega(x)^{-\theta} + 1), \quad \mbox{for} \ x \in \Omega, \ p \in \R^N,
\end{equation}
where $b_0, A > 0$ and $0\leq \theta<m$.

We first concentrate in regularity results in the {\bf superlinear} case
\begin{equation*}
m > \max \{ 1, \sigma \}, 
\end{equation*}
which encodes the coercivity of the Hamiltonian, see Theorems~\ref{teoholderbdy}
and~\ref{teoholderint} below. Note that, in Section~\ref{sublinearsubsec}, we state also a result in the sublinear
case and, in Section~\ref{sec:reglevyito}, we extend our results in the superlinear case to L\'evy-Ito operators.

Over the exponent $\theta$, we assume $0 \leq \theta < m$ in order to state the blow-up behavior 
at the boundary of the right-hand side. Thus, our arguments rely over the (more general) equation
\begin{equation}\label{eqgeneral}\tag{P'}
-I_x(u,x) + b_0 \ |Du(x)|^m =  A (d_\Omega(x)^{-\theta} + 1), \quad x \in \Omega. 
\end{equation}

In principle, due to the nonlocal nature of $I_x$, the function $u$ satisfying the above equation should be defined not only in $\Omega$ but
on the set
\begin{equation}\label{Omeganu}
\Omega_\nu =  \Omega \cup \bigcup_{x \in \bar{\Omega}} \{ x + \mathrm{supp} \{ \nu_x \} \},
\end{equation}
which, loosely speaking, represents the reachable set from $\Omega$ through $\nu$. 

The following result states the regularity up to the boundary for \textsl{subsolutions} of problem~\eqref{eqgeneral}.
\begin{teo}\label{teoholderbdy}
Let $\Omega \subseteq \R^N$ be a bounded domain, $A, b_0 > 0$ and $\sigma \in (0,2)$. Let $\{ \nu_x \}_{x \in \R^N}$ be a family of
measures satisfying (M1)-(M2) relative to $\sigma$, and $I_x$ defined as in~\eqref{operator} 
if $\sigma \geq 1$ and as ~\eqref{operatorsigmasmall} if $\sigma < 1$, associated to $\{ \nu_x \}_{x \in \R^N}$. 
Let $m > \max\{1, \sigma\}$, $\theta \in [0, m)$, and define
\begin{equation}\label{gammaholderbdy}
\gamma_0 = \min \{ (m - \sigma)/m, (m - \theta)/m \}. 
\end{equation}

Then, any bounded viscosity subsolution $u: \R^N \to \R$ to the problem~\eqref{eqgeneral}
is locally H\"older continuous in $\Omega$ with H\"older exponent $\gamma_0$ as in~\eqref{gammaholderbdy}, and 
H\"older seminorm depending on $\Omega$, the data 
and $\mathrm{osc}_{\Omega_\nu}(u)$, with $\Omega_\nu$ defined as in~\eqref{Omeganu}.

Moreover, if $\Omega$ has a $C^{1,1}$ boundary, then $u$ can be extended to $\bar{\Omega}$ as a H\"older continuous function of 
exponent $\gamma_0$.
\end{teo}

A second result states interior H\"older regularity for subsolutions of~\eqref{eqgeneral} with a H\"older exponent which is more natural 
to the balance between the order of the nonlocal operator and the Hamiltonian.
\begin{teo}\label{teoholderint}
Let $\Omega \subseteq \R^N$ be a bounded domain, $A, b_0 > 0$ and $\sigma \in (0,2)$. Let $\{ \nu_x \}_{x \in \R^N}$ be a family of
measures satisfying (M1)-(M2) relative to $\sigma$, and $I_x$ defined as in~\eqref{operator} 
if $\sigma \geq 1$ and as ~\eqref{operatorsigmasmall} if $\sigma < 1$, associated to $\{ \nu_x \}_{x \in \R^N}$. 
Let $m > \max \{1,  \sigma\}$ and $\theta \in [0, m)$. 
Define
\begin{equation}\label{gamma}
\begin{displaystyle}
\tilde{\gamma}_0 = \tilde{\gamma}_0(\sigma, m) = \left \{ \begin{array}{cl} (m - \sigma)/(m - 1) \quad & \mbox{if} \ \sigma > 1 \\
\in (0,1) \quad & \mbox{if} \ \sigma = 1 \\
1 \quad & \mbox{if} \ \sigma < 1,
\end{array} \right .
\end{displaystyle}
\end{equation}
and consider
\begin{equation}\label{gammaholderint}
\gamma_0 = \min \{ \tilde{\gamma}_0, (m - \theta)/m \}.
\end{equation}

Then, any bounded viscosity subsolution $u: \R^N \to \R$ to the equation~\eqref{eqgeneral}
is locally H\"older continuous in $\Omega$ with exponent $\gamma_0$ given by~\eqref{gammaholderint}, and 
H\"older seminorm depending on the data, $\Omega$ 
and $\mathrm{osc}_{\Omega_\nu}(u)$, where $\Omega_\nu$ is defined as in~\eqref{Omeganu}.
\end{teo}

Note that for the same data, $\gamma_0$ defined in~\eqref{gammaholderint} is always bigger or equal than $\gamma_0$ defined 
in~\eqref{gammaholderbdy}, and therefore, the interior H\"older exponent 
given by Theorem~\ref{teoholderint} is better
than the one given by Theorem~\ref{teoholderbdy}. 

\begin{remark}\label{rmkunbounded}
Theorems~\ref{teoholderbdy} and~\ref{teoholderint} can be extended to unbounded domains. In fact, if $\Omega$ is unbounded, 
arguing over a bounded set $\Omega' \subset \Omega$ we can apply the method used in the above theorems to conclude the corresponding 
local H\"older regularity results for $\Omega$. Moreover, if $\partial \Omega$ has uniform $C^{1,1}$ bounds, and if (M1)-(M2) hold 
with $C_R$ independent of $R$, then we have global H\"older estimates for 
bounded subsolutions to~\eqref{eqgeneral}, in the flavour of Theorem~\ref{teoholderbdy}.
\end{remark}

Since our aim is to include in our regularity results nonlocal operators of censored nature, 
we provide here a more accurate definition of such an operator. Recalling definition~\eqref{Omeganu}, 
we say that $I_x$ is of censored nature relative to $\Omega$ if the family $\{ \nu_x \}_{x \in \R^N}$ defining $I_x$ satisfies the condition
\begin{equation}\label{censoredcondition}
\Omega_\nu = \Omega. 
\end{equation}

The idea is to set up the problem to provide an unified proof of Theorem~\ref{teoholderbdy} for censored and noncensored operators.
This is possible after a ``censoring'' procedure we explain now.
Let $\{ \nu_x\}_{x \in \R^N}$ a family of L\'evy measures and $\Omega \subseteq \R^N$ an open set. For each 
$\xi \in \R^N$ we define the \textsl{censored measure respect to $\Omega$ and $\xi$} as
\begin{equation}\label{tildenu}
\tilde{\nu}_\xi(dz) = \mathbf{1}_{\Omega - \xi}(z) \nu_\xi(dz). 
\end{equation}

For $\xi, x \in \R^N$, $\delta > 0$ and a bounded function $\phi \in C^2(\bar{B}_\delta(x)),$ we define
\begin{equation}\label{tildeIcensored}
\begin{split}
\tilde{I}_\xi(\phi, x) = & \int_{\R^N} [\phi(x + z) - \phi(x) - \mathbf{1}_B \langle D\phi(x), z \rangle] \tilde{\nu}_\xi(dz) \\ 
= & \int_{\Omega - \xi} [\phi(x + z) - \phi(x) - \mathbf{1}_B \langle D\phi(x), z \rangle] \nu_\xi (dz).
\end{split}
\end{equation}

Of special interest is the \textsl{censored operator} 
$I_\Omega$ defined as
\begin{equation}\label{Icensored}
I_\Omega (\phi, x) = \tilde{I}_x(\phi, x), \quad x \in \bar{\Omega},
\end{equation}
from whose definition we note that $I_\Omega(\phi, x) = I_x[\Omega - x](\phi, x)$.

Note that if $\{ \nu_x \}_{x \in \R^N}$ satisfies (M1) and (M2), then $\{ \tilde{\nu}_x \}_{x\in \R^N}$ satisfies (M1) and (M2)
with the same constants $C_R$. Thus, the next lemma allows us to reduce general nonlocal equations like~\eqref{eqgeneral}
to the censored case.
\begin{lema}\label{lemacensored}{\bf (Censoring the Equation) }
Let $\Omega \subset \R^N$ open and bounded, $\sigma \in (0,2)$ and $\{ \nu_x \}_{x \in \R^N}$ a family of measures 
satisfying (M1)-(M2) related to $\sigma$. Let $I_x$ be as in~\eqref{operator},~\eqref{operatorsigmasmall} associated to $\{ \nu_x \}_{x \in \R^N}$.
Let $m > \sigma, \beta_0 > 0$ and 
for $f: \Omega \to \R$ locally bounded,  let $u : \R^N \to \R$ be a bounded viscosity subsolution to 
\begin{equation}\label{eqlemacensored}
-I_x(u,x) + \beta_0 |Du(x)|^m = f(x), \quad x \in \Omega.
\end{equation}

Then, there exists $C > 0$ (depending on $\Omega$ and $\beta_0$) such that the function $u$ restricted to $\Omega$ satisfies, in the viscosity sense, 
the inequality
\begin{equation*}\label{eqlemacensored2}
- I_\Omega(u,x) + \frac{\beta_0}{2} |Du(x)|^m \leq f(x) + C (\mathrm{osc}_{\Omega_\nu}(u) + 1)d_\Omega(x)^{-\sigma}, \quad x \in \Omega,
\end{equation*}
where $I_\Omega$ is defined in~\eqref{Icensored} and $\Omega_\nu$ is defined in~\eqref{Omeganu}.
\end{lema}

\noindent
{\bf \textit{Proof:}} For simplicity, we present the proof for classical subsolutions. The rigorous proof 
follows easily by using classical viscosity techniques
(for instance, see~\cite{Topp}).
We also focus on the case $\sigma \geq 1$.

Using~\eqref{eqlemacensored}, for each $x \in \Omega$ we have
\begin{equation*}
\begin{split}
& - I_\Omega(u,x) + \beta_0 |Du(x)|^m \\
\leq & \ f(x) + \int_{\Omega^c - x} (u(x + z) - u(x))\nu_x(dz) 
+ |Du(x)| \int_{B \cap (\Omega^c - x)} |z| \nu_x(dz) \\
\leq & \ f(x) + C \Big{(} \mathrm{osc}_{\Omega_\nu}(u) d_\Omega(x)^{-\sigma} + |Du(x)| h_{1, \sigma}(d_\Omega(x)) \Big{)},
\end{split}
\end{equation*}
where $C > 0$ comes from the application of (M1) and depends only on $\Omega$.
Now, by Young's inequality, there exists $C(\beta_0)$ such that
\begin{equation*}
|Du(x)| d_\Omega(x)^{1-\sigma} \leq \frac{\beta_0}{2} |Du(x)|^m + C(\beta_0) h_{1, \sigma}(d_\Omega(x))^{m/(m - 1)}.
\end{equation*}

At this point, we note that since $m > \sigma$ we have $m(1 - \sigma)/(m - 1) \geq - \sigma$.
Then, if $\sigma > 1$, using~\eqref{defh} we can write
\begin{equation*}
h_{1, \sigma}(d_\Omega(x))^{m/(m - 1)} = d_\Omega(x)^{m(1 - \sigma)/(m - 1)} \leq  d_\Omega(x)^{- \sigma}, 
\end{equation*}
meanwhile if $\sigma = 1$, we get
\begin{equation*}
h_{1, \sigma}(d_\Omega(x))^{m/(m - 1)} = (|\log(d_\Omega(x))| + 1)^{m/(m - 1)} \leq C d_\Omega(x)^{-\sigma},
\end{equation*}
where $C > 0$ depends only on $m$. Thus, using these estimates we conclude the result for the case $\sigma \geq 1$.

The case $\sigma < 1$ follows the same ideas but with easier computations because of the first order finite difference of the 
integrand defining $I_x$, see~\eqref{operatorsigmasmall}.
\qed


\subsection{Key Technical Lemmas.}
We start with some notation:
for $r > 0$ and $x_0 \in \R^N$, define
\begin{equation}\label{tau}
d_0(x) = |x - x_0| \quad \mbox{and} \quad d_r(x) = r - d_0(x),
\end{equation}
that is, for $x \in B_r(x_0)$, $d_0(x)$ represents the distance of $x$ to the center of the ball, meanwhile $d_r(x)  = d_{B_r(x_0)}(x)$ 
is the distance of $x$ to the boundary of the ball. We define $w$ as
\begin{equation}\label{defw}
w =  w_1 + w_2,
\end{equation}
where, for $C_1, \gamma > 0$ and $C_2 \geq 0$ we consider
\begin{equation}\label{defwi}
\begin{split}
w_1(x) & = \left \{ \begin{array}{ll} C_1 d_0(x)^\gamma \ & x \in \bar{B}_r(x_0) \\ 
C_1 r^\gamma \ & x \in \bar{B}_r^c(x_0) \end{array} \right . \\
w_2(x) & = \left \{ \begin{array}{ll} C_1 (r^\gamma - d_r(x)^\gamma) \ & x \in \bar{B}_r(x_0) \\ 
C_1 r^\gamma + C_2 \ & x \in \bar{B}_r^c(x_0). \end{array} \right . 
\end{split}
\end{equation}

We note that $w_1$ and $w_2$ (when $C_2 = 0$) are H\"older continuous  in $\R^N$ with exponent $\gamma$. If $C_2 > 0$, $w_2$ 
is $\gamma$-H\"older in $B_r(x_0)$ and it has a discontinuity on $\partial B_r(x_0)$. In any case, both $w_1$ and $w_2$ (for any $C_2 \geq 0$)
are smooth in $B_r(x_0) \setminus \{ 0 \}$.

For $x \in B_r(x_0)$ consider $\varrho$ defined as
\begin{equation}\label{varrho}
\varrho(x) = \frac{1}{4} \min \{ d_0(x), d_r(x) \}.
\end{equation}

Of course, $w$ depends on the particular choice of $\gamma, r, x_0, C_1, C_2$, meanwhile $\varrho$ depends on $r$ and $x_0$, but we omit 
these dependences for simplicity of the notation. 

We remark that if $|x - x_0| \leq r/2$ then $\varrho(x) = d_0(x)/4$, meanwhile if $|x - x_0| > r/2 $ we have $\varrho(x) = d_r(x) / 4$.

The goal is to prove that $w$ is a supersolution of~\eqref{eqgeneral}.
The following key lemma gives us a first useful estimate for the nonlocal term applied to $w$.
\begin{lema}\label{lemaIw}
Let $\sigma \in (0,2)$ and a family of measures $\{ \nu_x \}_{x \in \R^N}$ satisfying (M1), (M2) relative to $\sigma$. Let
$I_x$ as in~\eqref{operator},~\eqref{operatorsigmasmall} associated to $\{ \nu_x \}_{x \in \R^N}$.
Let $x_0 \in \R^N$, $r \in (0,1)$, $\gamma \in (0, 1]$, $C_1 > 0$, $C_2 \geq 0$, and consider $w$ as in~\eqref{defw} and $\varrho$ as in~\eqref{varrho}
associated to these parameters. Then, there exists a constant $C > 0$ (not depending on $r$, $C_1$ and $C_2$) such that
\begin{equation}\label{estimateIw}
\sup \limits_{ \xi \in B_1(x)} \{ I_\xi (w, x) \} 
\leq C \left \{ \begin{array}{ll} C_1  \varrho^{\gamma - 1}(x) h_{1, \sigma}(\varrho(x)) & \ \mbox{if} \ C_2 = 0, \sigma \geq 1 \\ 
C_1 h_{\gamma, \sigma}(\varrho(x)) & \ \mbox{if} \ C_2 = 0, \sigma < 1 \\
(C_1 + C_2) \varrho(x)^{- \sigma} & \ \mbox{if} \ C_2 > 0
\end{array} \right . ,
\end{equation} 
for each  $x \in B_r(x_0) \setminus \{ x_0 \}$. 
\end{lema}

\noindent
{\bf \textit{Proof:}}  Denote $R = |x_0| + 1$.
We remark that $C_R$ in the arguments to come is a generic constant depending on $R$ through the constants arising in (M1) and (M2). The constant 
$C$ arising in the proof is a positive constant independent of $x, R, r, C_1$ or $C_2$.

Consider $x \in B_r(x_0) \setminus \{ x_0 \}$. For each $\xi \in B_1(x)$, by definition of $w$ we can write 
$$
I_\xi(w, x) = I_\xi(w_1, x) + I_\xi(w_2, x),
$$   
where $w_i$, $i = 1,2$ are defined in~\eqref{defwi}. In what follows, we are going to estimate the integrals in the right-hand side of the above
expression.

\medskip
\noindent
1.- \textsl{ Estimate for $I_\xi(w_1, x)$.} We can split this integral term as
\begin{equation*}
I_\xi(w_1, x) = I_\xi[B_{\varrho(x)}](w_1, x) + I_\xi[B_{\varrho(x)}^c](w_1, x). 
\end{equation*}

Note that for each $z \in B_{\varrho(x)}$ we have
\begin{equation*}
\begin{split}
w_1(x + z) - w_1(x) = & \ \langle D w_1(x + t z), z \rangle, \\
w_1(x + z) - w_1(x) - \langle Dw_1(x), z\rangle = & \ \frac{1}{2}\langle D^2 w_1(x + s z) z, z \rangle, 
\end{split}
\end{equation*}
for some $s, t \in (0,1)$. We recall that the first equality is used in the integral defining 
$I_\xi[B_{\varrho(x)}](w_1, x)$ when $\sigma < 1$, and the 
second is used in the case $\sigma \geq 1$. Now, direct computations on the derivatives of $w_1$ drives us to 
\begin{equation*}
\begin{split}
\langle D^2 w_1(x + s z) z, z \rangle \leq & \ C_1 \gamma d_0(x)^{\gamma - 2} |z|^2 \\ 
\langle D w_1(x + t z), z \rangle \leq & \ C_1 \gamma d_0(x)^{\gamma - 1} |z|.
\end{split}
\end{equation*}
for all $z \in B_{\varrho(x)}$, $s, t \in (0,1)$. Thus, using these inequalities on the corresponding form of $I_\xi[B_{\varrho(x)}](w_1, x)$,
using that $\varrho(x) \leq d_0(x)$ and applying (M2), we arrive at
\begin{equation}\label{IBxw1}
I_\xi[B_{\varrho(x)}](w_1, x) \leq C_R C_1 \varrho(x)^{\gamma - \sigma}.
\end{equation}

Concerning the estimate of $I_\xi[B_{\varrho(x)}^c](w_1, x)$, we write
\begin{equation*}
\begin{split}
I_\xi[B_{\varrho(x)}^c](w_1, x) \leq \int \limits_{B_{\varrho(x)}^c} [w_1(x + z) - w_1(x)] \nu_\xi(dz) 
+ |D w_1(x)| \int \limits_{B \setminus B_{\varrho(x)}} |z| \nu_\xi(dz),
\end{split}
\end{equation*}
and we suppress the last integral term in the case $\sigma < 1$.
Using the definition of $w_1$ we get from the above inequality that
\begin{equation*}
\begin{split}
I_\xi[B_{\varrho(x)}^c](w_1, x) \leq & \ \int_{B \setminus B_{\varrho(x)}} [w_1(x + z) - w_1(x)] \nu_\xi(dz) 
+ C_1 r^\gamma \int_{B^c} \nu_\xi(dz) \\
& \ + C_1 \gamma d_0(x)^{\gamma - 1} \int_{B \setminus B_{\varrho(x)}} |z| \nu_\xi(dz),
\end{split}
\end{equation*}
where, as before, the last integral does not exist if $\sigma < 1$.
Since $w_1$ is $\gamma-$H\"older continuous we have $w_1(x + z) - w_1(x) \leq C_1 |z|^\gamma$. 
Using this together with (M1) (see~\eqref{defh}) we can write
\begin{equation*}
I_\xi[B_{\varrho(x)}^c](w_1, x) \leq C_R C_1 \Big{(} h_{\gamma, \sigma}(\varrho(x)) + r^\gamma + d_0(x)^{\gamma - 1} h_{1, \sigma}(\varrho(x)) \Big{)},
\end{equation*}
where the last term inside the parentheses is suppressed if $\sigma < 1$. Noting that $\varrho(x) \leq d_0(x) < r < 1$, we conclude that
\begin{equation*}
I_\xi[B_{\varrho(x)}^c](w_1, x) \leq C_R C_1 
\left \{ \begin{array}{ll} h_{\gamma, \sigma}(\varrho(x)) + \varrho(x)^{\gamma - 1} h_{1, \sigma}(\varrho(x)), \ & \mbox{if} \ \sigma \geq 1 \\
h_{\gamma, \sigma}(\varrho(x)), \ & \mbox{if} \ \sigma < 1. \end{array} \right .
\end{equation*}

At this point, we note that if $\sigma \geq 1$ and $\gamma \in (0, 1]$, we always have 
$h_{\gamma, \sigma}(\varrho) \leq \varrho^{\gamma - 1} h_{1, \sigma}(\varrho)$, for all $\varrho \in (0,1)$.
Taking this into account we get
\begin{equation*}
I_\xi[B_{\varrho(x)}^c](w_1, x) \leq C_R C_1 
\left \{ \begin{array}{ll} \varrho(x)^{\gamma - 1} h_{1, \sigma}(\varrho(x)), \ & \mbox{if} \ \sigma \geq 1 \\
h_{\gamma, \sigma}(\varrho(x)), \ & \mbox{if} \ \sigma < 1. \end{array} \right .
\end{equation*}
and joining this last inequality and~\eqref{IBxw1} we conclude that
\begin{equation}\label{Iw1}
I_\xi(w_1, x) \leq C_R C_1 
\left \{ \begin{array}{ll} \varrho(x)^{\gamma - 1} h_{1, \sigma}(\varrho(x)), \ & \mbox{if} \ \sigma \geq 1 \\
h_{\gamma, \sigma}(\varrho(x)), \ & \mbox{if} \ \sigma < 1. \end{array} \right .
\end{equation}


\medskip
\noindent
2.- \textsl{ Estimate for $I_\xi(w_2, x)$.} 
Analogously as the previous estimate, we write
\begin{equation*}
I_\xi(w_2, x) = I_\xi[B_{\varrho(x)}](w_2, x) + I_\xi[B_{\varrho(x)}^c](w_2, x). 
\end{equation*}

We start with $I_\xi[B_{\varrho(x)}](w_2, x)$. By recalling~\eqref{defwi}, direct computations drive us to
\begin{equation*}
\begin{split}
Dw_2(x) = & \ C_1 \gamma d_r^{\gamma - 1}(x) Dd_0(x), \\
D^2 w_2(x) 
= & \ C_1 \gamma d_r(x)^{\gamma - 2} d_0(x)^{-1} \\
& \ \times \Big{(} d_r(x) I_N + [(1 - \gamma)d_0(x) - d_r(x)] Dd_0(x) \otimes Dd_0(x) \Big{)},
\end{split}
\end{equation*}
and therefore, using the above computations as a Taylor expansion of the finite difference in the integral
defining $I_\xi[B_{\varrho(x)}](w_2, x)$, we claim that 
\begin{equation}\label{IBxw2}
I_\xi[B_{\varrho(x)}](w_2, x) \leq C_R C_1 \varrho(x)^{\gamma-\sigma}.
\end{equation}

In fact, when $\sigma < 1$, using~\eqref{operatorsigmasmall} and the above expression for $Dw_2$, we have
\begin{equation*}
I_\xi[B_{\varrho(x)}](w_2, x) = C_1 \gamma \int_{0}^{1} \int_{B_{\varrho(x)}} d_r^{\gamma - 1}(x + sz) \langle Dd_0(x + sz), z\rangle \nu_\xi(dz) ds, 
\end{equation*}
but for all $s \in (0,1)$ and $z \in B_{\varrho(x)}$, we have $d_r(x + sz) \geq \varrho(x)$. Thus, we have
\begin{equation*}
I_\xi[B_{\varrho(x)}](w_2, x) \leq C C_1 \varrho^{\gamma - 1}(x) \int_{B_{\varrho(x)}} |z| \nu_\xi(dz),  
\end{equation*}
and applying (M1) we conclude~\eqref{IBxw2}.

Now we deal with the case $\sigma \geq 1$. Since in this case
\begin{equation*}
\begin{split}
I_\xi[B_{\varrho(x)}](w_2, x) = \frac{1}{2} \int_{0}^{1} \int_{B_{\varrho(x)}} \langle D^2w_2(x + s z) z, z \rangle \nu(dz) ds,
\end{split}
\end{equation*}
using the explicit form of $D^2 w_2$ we get
\begin{equation}\label{IBxw2a}
\begin{split}
& I_\xi[B_{\varrho(x)}](w_2, x) \\
\leq & \ C C_1 r \int_{0}^{1} \int_{B_{\varrho(x)}} d_r(x + sz)^{\gamma - 2} d_0(x + sz)^{-1} |z|^2 \nu_\xi(dz) ds,
\end{split}
\end{equation}
and we estimate this last integral by cases. If $d_0(x) \geq r/2$ we have $\varrho(x) = d_r(x)/4$. Then, 
for $z \in B_{\varrho(x)}$ and $s \in (0,1)$ we have $3\varrho(x) \leq d_r(x + sz)$ and $r/4 \leq d_0(x + sz)$. 
Using these estimates into~\eqref{IBxw2a}, we conclude 
\begin{equation*}
\begin{split}
I_\xi[B_{\varrho(x)}](w_2, x)
\leq C C_1 \varrho(x)^{\gamma - 2} \int_{B_{\varrho(x)}} |z|^2 \nu_\xi(dz) \leq C_R C_1 \varrho(x)^{\gamma - \sigma},
\end{split}
\end{equation*}
where we have used (M2). On the other hand, if $d_0(x) < r/2$ we have $\varrho(x) = d_0(x)/4$. Then, for $z \in B_{\varrho(x)}$ and $s \in (0,1)$
we have $r/4 \leq d_r(x + sz)$ and $3\varrho(x) \leq d_0(x + sz)$. Using these estimates into~\eqref{IBxw2a}, we get
\begin{equation*}
I_\xi[B_{\varrho(x)}](w_2, x)
\leq C C_1 r^{\gamma - 1} \varrho(x)^{-1} \int_{B_{\varrho(x)}} |z|^2 \nu_\xi(dz) \leq C_R C_1 \varrho(x)^{\gamma - \sigma},
\end{equation*}
where we have used that $\varrho(x) \leq r$ and (M2). This concludes~\eqref{IBxw2}.

Concerning the estimate of $I_\xi[B_{\varrho(x)}^c](w_2, x)$, we should be careful with the fact that $C_2$ may be strictly positive.

At one hand, if $C_2 = 0$, then as in the computations relative to $w_1$, we have
\begin{equation*}
w_2(x + z) - w_2(x) \leq C_1 |z|^\gamma \quad \mbox{for all} \ z \in B_{\varrho(x)}^c,
\end{equation*}
and therefore, we can write
\begin{equation*}
\begin{split}
& I_\xi[B_{\varrho(x)}^c](w_2, x) \\
\leq & \ \int_{B_{\varrho(x)}^c} [w_2(x + z) - w_2(x)] \nu_\xi(dz) 
+ |D w_2(x)| \int_{B \setminus B_{\varrho(x)}} |z| \nu_\xi(dz) \\
\leq & \ C_1 \int_{B_{\varrho(x)}^c} |z|^\gamma \nu_\xi(dz) + C_1 d_r(x)^{\gamma - 1}\int_{B \setminus B_{\varrho(x)}} |z| \nu_\xi(dz),
\end{split}
\end{equation*}
where the last integral is suppressed if $\sigma < 1$.
Thus, applying (M1) and using that $r < 1$, we obtain from the above inequality that
\begin{equation*}
I_\xi[B_{\varrho(x)}^c](w_2, x) \leq C_R C_1 \Big{(} h_{\gamma, \sigma}(\varrho(x)) + d_r(x)^{\gamma - 1} h_{1, \sigma}(\varrho(x))\Big{)},
\end{equation*}
where the last term does not exist if $\sigma < 1$.
Finally, since $\varrho(x) \leq d_r(x)$ we conclude
\begin{equation}\label{IBxw2extC2=0}
I_\xi[B_{\varrho(x)}^c](w_2, x) 
\leq C_R \ C_1 \left \{ \begin{array}{ll} \varrho(x)^{\gamma - 1} h_{1, \sigma}(\varrho(x)), \ & \mbox{if} \ \sigma \geq 1 \\
h_{\gamma, \sigma}(\varrho(x)), \ & \mbox{if} \ \sigma < 1. \end{array} \right .
\end{equation}

On the other hand, if $C_2 > 0$, then we have the inequality
\begin{equation*}
w_2(x + z) - w_2(x) \leq C_1 + C_2 \quad \mbox{for all} \ z \in B_{\varrho(x)}^c.
\end{equation*}

Using this, now we can write
\begin{equation*}
\begin{split}
& I_\xi[B_{\varrho(x)}^c](w_2, x) \\
\leq & \ \int_{B_{\varrho(x)}^c} [w_2(x + z) - w_2(x)] \nu_\xi(dz) 
+ |D w_2(x)| \int_{B \setminus B_{\varrho(x)}} |z| \nu_\xi(dz) \\
\leq & \ (C_1 + C_2)\int_{B_{\varrho(x)}^c} \nu_\xi(dz) + C_1 d_r(x)^{\gamma - 1}\int_{B \setminus B_{\varrho(x)}} |z| \nu_\xi(dz),
\end{split}
\end{equation*}
where the last integral is suppressed if $\sigma < 1$.
Applying (M1) and using that $\varrho(x) \leq d_r(x)$ we conclude in this case that
\begin{equation*}\label{IBxw2extC2>0}
I_\xi[B_{\varrho(x)}^c](w_2, x) \leq C_R (C_1 + C_2) \varrho(x)^{- \sigma} + C_R C_1 \varrho^{\gamma - 1}(x) h_{1, \sigma}(\varrho(x)),
\end{equation*}
where the last term does not exist if $\sigma < 1$. Thus, since $\gamma > 0$ we get
\begin{equation}\label{IBxC2>0}
I_\xi[B_{\varrho(x)}^c](w_2, x) \leq C_R (C_1 + C_2) \varrho(x)^{- \sigma}. 
\end{equation}

In summary, when $C_2 = 0$, joining~\eqref{IBxw2extC2=0} and~\eqref{IBxw2} we have
\begin{equation}\label{Iw2C2=0}
I_\xi(w_2, x) \leq C_R \ C_1 \left \{ \begin{array}{ll} \varrho(x)^{\gamma - 1} h_{1, \sigma}(\varrho(x)), \ & \mbox{if} \ \sigma \geq 1 \\
h_{\gamma, \sigma}(\varrho(x)), \ & \mbox{if} \ \sigma < 1, \end{array} \right .
\end{equation}
meanwhile, when $C_2 > 0$, using~\eqref{IBxC2>0} and~\eqref{IBxw2} we conclude that 
\begin{equation}\label{Iw2C2>0}
I_\xi(w_2, x) \leq C_R (C_1 + C_2) \varrho(x)^{- \sigma}.
\end{equation}

\medskip
\noindent
3.- \textsl{Conclusion.} The estimate~\eqref{estimateIw} comes from~\eqref{Iw1} and~\eqref{Iw2C2=0} when $C_2 = 0$, and 
from~\eqref{Iw1} and~\eqref{Iw2C2>0} when $C_2 > 0$. The proof is complete.
\qed


Using the last lemma we are able to prove $w$ is a strict supersolution for a 
problem ad-hoc to~\eqref{eqgeneral}. 
This is established in the following two lemmas, whose main difference is whether $C_2$ is strictly positive or not.
\begin{lema}\label{keylemmabdy}{\bf ($\mathrm{\mbox{Strict Supersolution, Case} \ \mathbf{C_2 > 0}}$) }
Let $x_0 \in \R^N$, $\sigma \in (0,2)$ and $\{ \nu_x\}_{x \in \R^N}$ a family of measures satisfying (M1), (M2) relative to $\sigma$.
Consider $I_x$ as in~\eqref{operator},~\eqref{operatorsigmasmall} associated to $\{ \nu_x\}_{x \in \R^N}$.
Let $m > \max \{1, \sigma \}$, $\theta \in [0, m)$
and $\gamma_0$ given in~\eqref{gammaholderbdy}.

Then, for each $A, b_0, C_2 > 0$, there exists $C_1 > 0$ large enough such that, 
for all $r \in (0,1)$ and $\gamma \in (0, \gamma_0]$, the function $w$ defined in~\eqref{defw} (relative to $x_0, \gamma, C_1, C_2$ and $r$) 
satisfies the inequality
\begin{equation}\label{eqkeylemmabdy}
- \sup \limits_{\xi \in B_1(x)} \{ I_\xi (w, x) \} + b_0 |Dw(x)|^m \geq A \varrho(x)^{-\theta}
\quad \mbox{for} \ x \in B_r(x_0) \setminus \{ x_0 \}, 
\end{equation}
where $\varrho$ defined in~\eqref{varrho} is associated to $x_0$ and $r$.
\end{lema}

\noindent
{\bf \textit{Proof:}} Let $x \in B_r(x_0) \setminus \{ x_0 \}$. Direct computations over $w_1, w_2$ defined in~\eqref{defwi} 
give us the expression
\begin{equation*}
Dw(x) = C_1 \gamma (d_0(x)^{\gamma - 1} + d_r(x)^{\gamma - 1}) \frac{x - x_0}{|x - x_0|},
\end{equation*}
concluding that
\begin{equation*}\label{Domega}
|D w(x)| = C_1 \gamma (d_0(x)^{\gamma - 1} + d_r(x)^{\gamma - 1}) \geq C C_1 \varrho(x)^{\gamma - 1}.
\end{equation*}

Using this together with the estimates given by Lemma~\ref{lemaIw} for the nonlocal term in the case $C_2 > 0$,
we obtain the existence of an universal constant $\bar{C} > 0$ such that for all $C_1, C_2$ and $b_0$, and for all $x \in B_r(x_0) \setminus \{ x_0 \}$ 
we have
\begin{equation}\label{eqwkeylemmabdy}
\begin{split}
& -\sup \limits_{\xi \in B_1(x)} \{ I_\xi(w, x) \} + b_0|Dw(x)|^m \\
\geq & \ \bar{C} \Big{(} b_0 C_1^m \varrho(x)^{m(\gamma - 1)} - (C_1 + C_2) C_R \varrho(x)^{-\sigma} \Big{)}.
\end{split}
\end{equation}

But since $\gamma_0 = \min \{ m - \sigma, m - \theta \} /m$ and $\gamma \leq \gamma_0$ we have $m(\gamma - 1) \leq \min \{-\sigma, -\theta\}$. 
Then, we conclude from~\eqref{eqwkeylemmabdy} that
\begin{equation*}
-\sup \limits_{\xi \in B_1(x)} \{ I_\xi(w, x) \} + b_0|Dw(x)|^m 
\geq \bar{C} \varrho(x)^{m(\gamma - 1)} \Big{(} b_0 C_1^m - (C_1 + C_2) C_R \Big{)}.
\end{equation*}

Hence, we arrive at~\eqref{eqkeylemmabdy} by taking 
\begin{equation*}
C_1 = (4A(\bar{C}b_0)^{-1})^{1/m} + (4C_2 C_R b_0^{-1})^{1/m} + (2C_R b_0^{-1})^{1/(m - 1)},
\end{equation*}
that is, we should take $C_1$ satisfying
\begin{equation}\label{C1choicekeylemma2}
C_1 \geq C (A^{1/m} + C_2^{1/m} + 1), 
\end{equation}
where $C > 0$ is a constant not depending on $C_2$ or $A$.
\qed

Next lemma deals with the case $C_2 = 0$. 
\begin{lema}\label{keylemmaint}{\bf ($\mathrm{\mbox{Strict Supersolution, Case} \ \mathbf{C_2 = 0}}$) }
Let $x_0 \in \R^N$, $\sigma \in (0,2)$ and $\{ \nu_x\}_{x \in \R^N}$ a family of measures satisfying (M1), (M2) relative to $\sigma$.
Consider $I_x$ as in~\eqref{operator},~\eqref{operatorsigmasmall} associated to $\{ \nu_x\}_{x \in \R^N}$.
 Let $m > \max \{1, \sigma \}$, $\theta \in [0, m)$
and $\gamma_0$ defined in~\eqref{gammaholderint}.
Assume $C_2 = 0$. 
 
Then, for each $A, b_0 > 0$, there exists $C_1 > 0$ large enough such that, 
for all $r \in (0,1)$ and $\gamma \in (0, \gamma_0]$, the function $w$ defined in~\eqref{defw} (relative to $x_0, \gamma, C_1$ and $r$) 
satisfies the inequality~\eqref{eqkeylemmabdy}.
\end{lema}

The proof of this lemma follows exactly as Lemma~\ref{keylemmabdy} using the estimate given by Lemma~\ref{lemaIw} in the case $C_2 = 0$ 
and the definition of $\gamma_0$ given in~\eqref{gammaholderint}.

\begin{remark}\label{rmkgammas}
As we mentioned in the introduction, the power profile of $w$ gives us the H\"older regularity for subsolutions to~\eqref{eqgeneral}. The 
different uses of Lemmas~\ref{keylemmabdy} and~\ref{keylemmaint} can be described as follows:
as it can be seen in the proof of Theorem~\ref{teoholderbdy} below, the application of Lemma~\ref{keylemmabdy} under a correct choice of 
$C_2 > 0$ allows us to localize the arguments to obtain an interior H\"older regularity with a H\"older seminorm (cast by $C_1$) which 
is independent of the distance to the boundary, a key fact to conclude the regularity up to the boundary. However, the 
discontinuity of $w$ due to $C_2 > 0$
implies a ``worse'' bound for $I_x(w)$ (see Lemma~\ref{lemaIw}), restricting the values of the H\"older exponent if we look for
regularity up to the boundary, no matter the nonlocal operator has censored nature or not. 

On the other hand, Lemma~\ref{keylemmaint}
is used in the proof of Theorem~\ref{teoholderint}, where no localization is needed. Thus, the ``better'' bounds for $I_x(w)$
given by Lemma~\ref{lemaIw} allows to obtain interior H\"older regularity with ``more natural'' exponents.
\end{remark}


\subsection{Proofs of the Main Theorems.}

We start with the regularity result up to the boundary.

\medskip
\noindent
{\bf \textit{Proof of Theorem~\ref{teoholderbdy}:}} 
Applying Lemma~\ref{lemacensored}, we see that $u$ satisfies the censored equation
\begin{equation*}
- I_\Omega(u,x) + \frac{b_0}{2} |Du|^m \leq A (d_\Omega(x)^{-\theta} + 1) 
+ C(\mathrm{osc}_{\Omega_\nu}(u) + 1) d_\Omega(x)^{-\sigma}, \quad x \in \Omega,
\end{equation*}
where $C > 0$ is the constant given in Lemma~\ref{lemacensored}. If we define
$\eta = \max \{ \sigma, \theta \}$, in particular we see that $u$ satisfies the viscosity inequality
\begin{equation}\label{eqauxteoholderbdy}
- I_\Omega(u,x) + \frac{b_0}{2} |Du|^m \leq \tilde{A} d_\Omega(x)^{-\eta}, \quad x \in \Omega,
\end{equation}
where 
\begin{equation}\label{C1holderbdy1}
\tilde{A} = A(1 + \mathrm{diam}(\Omega)^\eta) + C(\mathrm{osc}_{\Omega_\nu}(u) + 1).
\end{equation}

From this point, we will argue over equation~\eqref{eqauxteoholderbdy}.

Let $x_0 \in \Omega$ and denote $R = |x_0| + 1$. Consider $\gamma_0$ as in~\eqref{gammaholderbdy}, and for $C_1, C_2 > 0$ to be fixed
later and $r = \min \{1,  d_\Omega(x_0) \}/4$, consider $w$ as in~\eqref{defw} (with $\gamma = \gamma_0$) associated to these parameters. 

Denote
\begin{equation*}\label{Mholderbdy}
M := \sup  \{ u(x) - u(x_0) - w(x) : x \in \bar{\Omega} \}. 
\end{equation*}

The aim is to prove that for suitable $C_1 > 0$ we get $M \leq 0$, which implies easily the H\"older continuity of $u$. 
We argue by contradiction, assuming that $M > 0$.
Choosing 
\begin{equation}\label{C2holderbdy}
C_2 \geq \mathrm{osc}_{\Omega_\nu}(u), 
\end{equation}
by definition of $w$, for each $x \in \bar{\Omega} \setminus \bar{B}_r(x_0)$ we have 
\begin{equation*}
u(x) - u(x_0) - w(x) \leq \mathrm{osc}_{\Omega_\nu}(u) - (2C_1 r^\gamma + \mathrm{osc}_{\Omega_\nu}(u)) < 0.
\end{equation*}

Hence, by the upper semicontinuity of $u -w$, it follows that
the supremum defining $M$ is attained in $\bar{B}_r(x_0)$. Moreover, since $w(x_0) = 0$, the point attaining the maximum in $M$
is in $\bar{B}_r(x_0) \setminus \{ x_0 \}$.

Let $A_0 > 0$ be fixed later. By Lemma~\ref{keylemmabdy}, we can consider $C_1$ large enough in order to have
\begin{equation}\label{eqteobdyw}
\begin{split}
- \sup \limits_{\xi \in B_1(x)} \{ \tilde{I}_\xi (w, x) \} + \frac{b_0}{2} |D w(x)|^m 
\geq  A_0 \varrho(x)^{-\eta},
\quad x \in B_{r}(x_0) \setminus \{ x_0 \},
\end{split}
\end{equation}
in fact, by~\eqref{C1choicekeylemma2} it is sufficient to take 
\begin{equation}\label{C1holderbdy2}
C_1 \geq C(A_0^{1/m} + C_2^{1/m}  + 1) 
\end{equation}
for some universal constant $C > 0$.
Doubling variables and penalizing, we consider
\begin{equation*}
M_\epsilon := \sup  \{ \Phi(x,y) : (x, y) \in \bar{\Omega} \times \bar{\Omega} \},
\end{equation*}
where $\Phi(x,y) = u(x) - u(x_0) - w(y) - \epsilon^{-2}|x - y|^2$. 

By classical arguments in the viscosity theory, we have $M_\epsilon \geq M > 0$ for all $\epsilon > 0$ and
the supremum in $M_\epsilon$ is attained at $(\bar{x}, \bar{y}) \in \bar{\Omega} \times \bar{\Omega}$ with 
$\bar{y} \in \bar{B}_r(x_0) \setminus \{ x_0 \}$, which in addition satisfies the following properties
\begin{equation}\label{class-visc}
\epsilon^{-2}|\bar{x}- \bar{y}|^2 \to 0; \quad \bar{x}, \bar{y} \to x^*; \quad u(\bar{x}) \to u(x^*), \quad \mbox{as} \ \epsilon \to 0,
\end{equation}
where $x^* \in \bar{B}_r(x_0) \setminus \{x_0 \}$ attains the supremum defining $M$. 
In particular, $\bar{y} \neq x_0$ for all $\epsilon > 0$. Moreover, note that the function
\begin{equation*}
-\Phi(\bar{x},\cdot ) :\, 
y \mapsto w(y) - (u(\bar{x}) - u(x_0) - \epsilon^{-2}|\bar{x} - y|^2) 
\end{equation*}
has a global minimum point at $\bar{y} \in \bar{B}_r(x_0) \setminus \{ x_0 \}$ for all $\epsilon > 0$. We claim that this fact implies
$\bar{y} \notin \partial B_r(x_0)$ for each $\epsilon > 0$. Otherwise, denoting $\xi = (x_0 - \bar{y})/|x_0 - \bar{y}|$ we have
$\bar{y} + s\xi \in B_r(x_0)$  for each $0 < s < r$. Therefore
$-\Phi(\bar{x},\bar{y})\leq -\Phi(\bar{x},\bar{y}+s\xi),$ which implies by definition of $w$ in~\eqref{defw}
\begin{equation*}
0 \leq s^{-1} (w(\bar{y}) - w(\bar{y} + s\xi))  \leq \epsilon^{-2} (-2\langle \bar{x} - \bar{y}, \xi \rangle + s)
\end{equation*}
and
\begin{equation*}
0 \leq  C_1  (s^{-1} (r^{\gamma} - (r - s)^\gamma) + s^{\gamma - 1}) 
\leq \epsilon^{-2} (-2\langle \bar{x} - \bar{y}, \xi \rangle + s).
\end{equation*}

Making $s \to 0$ we arrive at a contradiction, concluding the claim.
Hence, for all $\epsilon > 0$, there exists $r_\epsilon \in (0, r)$ such that $r_\epsilon < |\bar{y} - x_0| < r - r_\epsilon$.

On the other hand, using that $(\bar{x}, \bar{y})$ is a maximum point for $\Phi$, denoting $h = x - y$ and $\bar{h} = \bar{x} - \bar{y}$ 
we have 
\begin{equation*}
u(h + y) - w(y) - \epsilon^{-2}|h|^2 \leq u(\bar{h} + \bar{y}) - w(\bar{y}) - \epsilon^{-2}|\bar{h}|^2,
\end{equation*}
for each $y \in \bar{\Omega}$ and $h$ such that $y + h \in \bar{\Omega}$. Hence, we conclude 
\begin{equation*}
\bar{u}(y) - w(y) \leq \bar{u}(\bar{y}) - w(\bar{y}) \quad \mbox{for all} \ y \in \Omega - \bar{h},
\end{equation*}
where $\bar{u}(y) := u(\bar{h} + y)$ for each $y \in \Omega - \bar{h}$. In particular, $\bar{y}$ is a maximum point for $\bar{u} -w$ 
in $\Omega - \bar{h}$. Now, a simple translation argument over equation~\eqref{eqauxteoholderbdy} allows us to prove that $\bar{u}$ satisfies 
the equation
\begin{equation*}
- \tilde{I}_{x + \bar{h}}(\bar{u}, x) + \frac{b_0}{2} |D \bar{u}(x)|^m 
\leq \tilde{A} d_\Omega^{-\eta}(x + \bar{h}), \quad x \in \Omega - \bar{h},
\end{equation*}
in the viscosity sense.
Since $|\bar{h}| \to 0$ as $\epsilon \to 0$, for all $\epsilon$ small enough we have $\bar{y} \in B_r(x_0) \subset \Omega - \bar{h}$. 
Recalling $w$ is smooth at $\bar{y}$ we can use it as a test function for $\bar{u}$ at $\bar{y}$, concluding the inequality
\begin{equation*}
- \tilde{I}_{\bar{y} + \bar{h}}(w, \bar{y}) + \frac{b_0}{2} |D w (\bar{y})|^m \leq \tilde{A} d_\Omega^{-\eta}(\bar{y} + \bar{h}), 
\end{equation*}
but since $\bar{y} + \bar{h} \in B_1(\bar{y})$ for $\epsilon$ small enough, using~\eqref{eqteobdyw} we get
\begin{equation*}
A_0\varrho^{-\eta}(\bar{y}) \leq \tilde{A} d_\Omega^{-\eta}(\bar{y} + \bar{h}).
\end{equation*}

Note that for each $x \in B_r(x_0)$ we have $\varrho(x) \leq d_\Omega(x)$ and since $\eta \geq 0$, we get from the above inequality that
\begin{equation*}
A_0 d_\Omega^{-\eta}(\bar{y}) \leq \tilde{A}d_\Omega^{-\eta}(\bar{y} + \bar{h}).
\end{equation*}

At this point, recalling $\bar{h} \to 0$ and $\bar{y} \to x^* \in \bar{B}_r(x_0)$ as $\epsilon \to 0$, taking limits in the above inequality we 
arrive at a contradiction previously fixing 
\begin{equation}\label{C1holderbdy3}
A_0 \geq \tilde{A} + 1. 
\end{equation}

Thus, for each $x_0 \in \Omega$ and $r \leq d_\Omega(x_0)/4$, we have
\begin{equation*}\label{localHolderineq}
|u(x) - u(y)| \leq C_1 |x - y|^{\gamma_0} \quad \mbox{for all} \ x, y \in B_r(x_0),
\end{equation*}
from which we conclude the local H\"older continuity. In the case the boundary is $C^{1,1}$, from the above inequality we note that for each 
$B_r(x_0) \subset \Omega$, the H\"older exponent and seminorm of $u$ in $B_r(x_0)$ does not depend on $r$, and applying the method used by Barles
in~\cite{Barles1} (see also~\cite{Capuzzo-Dolcetta-Leoni-Porretta}) we can extend the H\"older regularity up to the boundary.

Finally, we recall that by~\eqref{C1holderbdy3},~\eqref{C1holderbdy1},~\eqref{C1holderbdy2} and the choice 
of $C_2$ in~\eqref{C2holderbdy}, the constant $C_1$ leading to the contradiction has the form
\begin{equation}\label{C1holderbdyfinal}
C_1 \geq C(A^{1/m} + \mathrm{osc}_{\Omega_\nu}(u)^{1/m} + 1), 
\end{equation}
for some constant $C > 0$ depending on the data. 
\qed

A very important consequence of the previous result is the following control of the oscillation.
\begin{cor}\label{corosc}{\bf (Oscillation Bound)}
Let $\Omega \subset \R^N$ be open and bounded with a $C^{1,1}$ boundary,
 and assume the hypotheses of Theorem~\ref{teoholderbdy} hold. Assume further the nonlocal operator has a 
censored nature, that is, the family of measures $\{ \nu_x \}_{x \in \R^N}$ satisfies the censored condition~\eqref{censoredcondition}. 
Then, there exists $K > 0$ such that, for each bounded viscosity subsolution $u$  of~\eqref{eqgeneral}, we have
\begin{equation*}
\mathrm{osc}_\Omega(u) \leq K. 
\end{equation*}
\end{cor}

\noindent
{\bf \textit{Proof:}} The choice of $C_1$ given by~\eqref{C1holderbdyfinal} in Theorem~\ref{teoholderbdy} leads us to
\begin{equation*}
|u(x) - u(y)| \leq C ( A^{1/m} + \mathrm{osc}_{\Omega_\nu}(u)^{1/m} + 1) |x - y|^{\gamma_0}, \quad \mbox{for all} \ x, y \in \bar{\Omega},
\end{equation*}
where $\gamma_0$ is given by~\eqref{gammaholderbdy}.
Now, by~\eqref{censoredcondition} we have $\mathrm{osc}_{\Omega_\nu}(u) = \mathrm{osc}_{\Omega}(u)$ and 
by compactness of $\bar{\Omega},$
there exists $\underline{x}, \bar{x} \in \bar{\Omega}$ such 
that $\mathrm{osc}_{\Omega}(u) = u(\bar{x})- u(\underline{x})$. Then, we can write
\begin{equation*}
\mathrm{osc}_{\Omega}(u) \leq C (A^{1/m} + \mathrm{osc}_\Omega(u)^{1/m} + 1),
\end{equation*}
from where we obtain the result since $m > 1$. 
\qed

Note that for noncensored problems, we can provide global oscillation bounds as in the last corollary if we a priori know that 
$\mathrm{osc}_{\Omega_\nu}(u) = \mathrm{osc}_{\Omega}(u)$.


\medskip
\noindent
{\bf \textit{Proof of Theorem~\ref{teoholderint}:}} 
Let $x_0 \in \Omega$, denote $R = |x_0| + 1$ and fix $r = \min \{ 1, d_\Omega(x_0) \}/4 $. 
Consider $\gamma_0$ as in~\eqref{gammaholderint} and for $C_1 > 0$ to be fixed later, define $w$ as in~\eqref{defw} (with $\gamma = \gamma_0$)
associated to these parameters.

Since the proof follows the same lines of Theorem~\ref{teoholderbdy}, we will be sketchy in the current proof bringing light on its contrasts.
The first difference is that this time we do not censorize the equation (since it would restrict the H\"older exponent, see Lemma~\ref{lemacensored}). 

Denote
\begin{equation}\label{Mteoint}
M := \sup  \{ u(x) - u(x_0) - w(x) : x \in \R^N \}. 
\end{equation}

The aim is to prove that for suitable $C_1 > 0$ we get $M \leq 0$. We argue by contradiction, assuming that $M > 0$.
Note that choosing 
\begin{equation}\label{C1oscteoint}
C_1 r^{\gamma_0} \geq \mathrm{osc}_{\Omega_\nu}(u), 
\end{equation}
and by the upper semicontinuity of $u - w$ we have the supremum defining $M$ is attained in $\bar{B}_r(x_0)$.

Let $A_0 > 0$ be fixed later. Enlarging $C_1$ if it is necessary, by Lemma~\ref{keylemmaint} we can write
\begin{equation}\label{eqteointw}
\begin{split}
- \sup \limits_{\xi \in B_1(x)} \{ I_\xi (w, x) \} + b_0 |D w(x)|^m 
\geq  A_0 \varrho(x)^{-\theta},
\quad x \in B_{r}(x_0) \setminus \{ x_0 \}.
\end{split}
\end{equation}
 
Doubling variables and penalizing, we consider
\begin{equation*}
M_\epsilon := \sup  \{ \Phi(x,y) : (x, y) \in \R^N \times \R^N \},
\end{equation*}
where $\Phi(x,y) = u(x) - u(x_0) - w(y) - \epsilon^{-2}|x - y|^2$. 
By classical arguments in the viscosity theory, we have $M_\epsilon \geq M > 0$ for all $\epsilon > 0$ and
the supremum in $M_\epsilon$ is attained at $(\bar{x}, \bar{y})$ with $\bar{x},\bar{y} \in \R^N$ with $\bar{y} \in \bar{B}_r(x_0) \setminus \{ x_0 \}$, 
which in addition satisfies~\eqref{class-visc}
where $x^* \in B_r(x_0) \setminus \{x_0 \}$ attains the supremum in~\eqref{Mteoint}.

If $\gamma_0 < 1$, then we can prove that $\bar{y} \notin \partial B_r(x_0)$ in the same way as in Theorem~\ref{teoholderbdy} using that
$w$ satisfies a state constraint problem on $\partial B_r(x_0)$. If $\gamma_0 = 1$ (which is the case of $\theta = 0$ and $\sigma < 1$), then 
we consider $w$ with $\gamma < \gamma_0$ and continue with the proof, taking into account that the H\"older seminorm does not change 
as $\gamma \to \gamma_0$.

From this point, we follow the remaining lines of Theorem~\ref{teoholderbdy}, taking $A_0$ large in terms of $A$ arising in~\eqref{eqgeneral} .
\qed


\subsection{Examples.}
\label{examplesubsection} 
In this section we provide some examples of nonlocal terms and Hamiltonians for which our results hold.

We start with the assumptions over the nonlocal term. As we mentioned before, assumptions (M1) and (M2) are intended as a restriction on the 
order of the operator, which is less or equal than $\sigma$. In the case of $x$-independent operators, that is the case when there exists a 
measure $\nu$ such that the family $\{ \nu_x \}_{x \in \R^N}$ defining $I_x$ satisfies $\nu_x = \nu$ for each $x \in \R^N$,
the operator may range from zero order operators (when $\nu$ is finite, see~\cite{Chasseigne}) to the factional Laplacian of order $s$ 
for $s \leq \sigma$, passing through operators which are not uniformly elliptic in the sense of Caffarelli and Silvestre~\cite{Caffarelli-Silvestre}, 
as it is the case of measures with the form
\begin{equation*}
\nu(dz) =  \mathbf{1}_{\mathbb{H}_+}(z) |z|^{-(N + s)} dz,
\end{equation*}
where $0 < s \leq \sigma$ and $\mathbb{H}_+ = \{ (z', z_N) \in \R^N : z_N > 0 \}$. Another interesting example of such non-uniformly 
elliptic operators is given by operators with ``orthogonal diffusion'', for example in the case $\nu$ has the form
\begin{equation}\label{crossedmeasure}
\nu(dz) = |z_2|^{-(N + s_2)} dz_2 \otimes \delta_{0}(z_1) dz_1 + |z_1|^{-(N + s_1)} dz_1 \otimes \delta_{0}(z_2) dz_2
\end{equation}
where $z = (z_1, z_2)$ with $z_i \in \R^{d_i}, i=1,2$ and $N = d_1 + d_2$, and $0 < s_1, s_2 \leq \sigma$. Here $\delta_0$ denotes the Dirac measure
supported at $0$ and $\otimes$ denotes the measure product. In this case, such a measure gives rise to an operator which is the sum of fractional 
Laplacians in each direction $z_i, i=1,2$.

Concerning $x$-dependent nonlocal operators, the classical example comes from measures $\nu_x$ with the form
\begin{equation*}
\nu_x(dz) = K(x, z) \nu(dz), 
\end{equation*}
where $\nu$ is an $x$-independent L\'evy measure and $K : \R^N \times \R^N \to \R$ is a nonnegative function such
that $K(\cdot, z) \in L^\infty_{loc}(\R^N)$
for all $x \in \R^N$, and $K(x, \cdot) \in L^\infty(\R^N)$ for all $x \in \R^N$. As a particular case we have the weighted fractional Laplacian
\begin{equation*}
-I_x(u,x) = K(x) (-\Delta)^\sigma u(x),
\end{equation*}
where $K$ is bounded and nonnegative.

We highlight that in view of Lemma~\ref{lemacensored}, the regularity results apply to censored operators defined in~\eqref{Icensored}, where 
we recall that the measures defining them has the form~\eqref{tildenu}.

Concerning $H$, we note that the structure of the Hamiltonian is encoded by 
the inequality~\eqref{H}. Thus, given $\sigma$ and $m > \max \{ 1, \sigma \}$, our results apply to $H$ with the form
\begin{equation}\label{Hexample1}
H(x,p) = b(x) |p|^m + a_1(x) |p|^l + \langle a_2(x), p \rangle - f(x),
\end{equation}
where $x, p \in \R^N$, $b \geq b_0 > 0$, $0 < l < m$ and $a_1, a_2, f$ bounded. In the case $m \leq 1$ we can consider
\begin{equation}\label{Hexample2}
H(x,p) = b(x) |p|^m + a_1(x) |p|^l - f(x),
\end{equation}
with $b, a_1, l$ and $f$ as above.
Of course, we can replace the main power $|p|^m$ by $\phi(x, p) |p|^m$, where the function $\phi : \R^N \times \R^N \to \R$ 
satisfies $\phi \geq \phi_0$ for some constant $\phi_0 > 0$.


\subsection{Regularity Results for the Sublinear Case.}
\label{sublinearsubsec}
In this subsection we provide a regularity results in the case $\sigma < m \leq 1$. 
\begin{teo}\label{teosublinearbdy}
Let $\Omega \subseteq \R^N$ be a bounded domain and $\sigma \in (0,1)$. Let $I_x$ as in~\eqref{operatorsigmasmall} associated to a family of measures 
$\{\nu_x\}_{x \in \R^N}$ satisfying (M1), (M2) relative to $\sigma$. Let $m \in (\sigma, 1]$, $\theta \in [0,m)$ and 
$\gamma_0$ as in~\eqref{gammaholderbdy}.

Then, for each $b_0, A > 0$ and $\gamma < \gamma_0$, any bounded viscosity subsolution $u: \R^N \to \R$ to the equation~\eqref{eqgeneral}
is locally H\"older continuous in $\Omega$ with H\"older exponent $\gamma$. 
If $\Omega$ has $C^{1,1}$ boundary, then $u$ is $\gamma$-H\"older continuous in $\Omega$ and can be extended as a H\"older continuous 
function on $\bar{\Omega}$.

The H\"older seminorm depends on the data and $\mathrm{osc}_{\Omega_\nu}(u)$, where $\Omega_\nu$ is defined in~\eqref{Omeganu}.
\end{teo}

\noindent
{\bf \textit{Proof:}} As in Theorem~\ref{teoholderbdy}, we start with the analogous of Lemma~\ref{keylemmabdy}. Let $r > 0$, consider
$x_0 \in \R^N$, define $d_0, d_r$ as in~\eqref{tau} and $\varrho$ as in~\eqref{varrho}. Let $w$ defined in~\eqref{defw}
associated to these parameters and $\gamma < \gamma_0$. Let $A, b_0 > 0$. Performing the same computations as in Lemma~\ref{keylemmabdy}
we arive at inequality~\eqref{eqwkeylemmabdy}, that is
\begin{equation*}
\begin{split}
& -\sup \limits_{\xi \in B_1(x)} \{ I_\xi(w, x) \} + b_0|Dw(x)|^m \\
\geq & \ \bar{C} \Big{(} b_0 C_1^m \varrho(x)^{m(\gamma - 1)} - (C_1 + C_2) C_R \varrho(x)^{-\sigma} \Big{)},
\end{split}
\end{equation*}
for all $x \in B_r(x_0) \setminus \{ x_0 \}$. Since this time $m(\gamma - 1) < -\sigma$ and $\varrho(x) \leq r$ for each $x \in B_r(x_0)$, 
we can take $r=r(C_1, C_2, b_0)$ small such that
\begin{equation*}
-\sup \limits_{\xi \in B_1(x)} \{ I_\xi(w, x) \} + b_0|Dw(x)|^m 
\geq \frac{\bar{C} b_0 C_1^m}{2} \varrho(x)^{m(\gamma - 1)}, \quad x \in B_r(x_0) \setminus \{ x_0 \}.
\end{equation*}

By the choice of $\gamma < \gamma_0$, we see that $m(\gamma -1) \leq -\theta$, and therefore 
$w$ satisfies
\begin{equation*}
-\sup \limits_{\xi \in B_1(x)} \{ I_\xi(w, x) \} + b_0|Dw(x)|^m 
\geq \bar{C} C_1^m \varrho(x)^{-\eta}, \quad x \in B_r(x_0) \setminus \{ x_0 \},
\end{equation*}
with $\eta = \max \{ \sigma, \theta \}$. From this point, we proceed exactly as in the proof of Theorem~\ref{teoholderbdy}, 
where the last inequality plays the role of~\eqref{eqteobdyw}, concluding the result by taking $C_1$ large in terms of $A$. 
\qed

\begin{remark}\label{rmksublinear}
Since $m \leq 1$, the parameter $r$ depends on $C_2$ in the proof of Theorem~\ref{teosublinearbdy} and therefore we have a H\"older seminorm
which does not give a control of the oscillation in the general case.
\end{remark}

Interior regularity results for the sublinear case in the flavour of Theorem~\ref{teoholderint} can be obtained in the same way as the previous 
theorem.
\begin{teo}\label{teosublinearint}
Let $\Omega \subseteq \R^N$ be a bounded domain.
Let $\sigma \in (0,1)$, $I_x$ as in~\eqref{operatorsigmasmall} associated to a family of measures 
$\{ \nu_x \}_{x \in \R^N}$ satisfying (M1), (M2) relative to $\sigma$. Let $m \in (\sigma, 1]$, $\theta \in [0, m)$ and 
$\gamma_0$ as in~\eqref{gammaholderint}.

Then, for each $b_0, A > 0$ and $\gamma < \gamma_0$, any bounded viscosity subsolution $u: \R^N \to \R$ to the equation~\eqref{eqgeneral}
is locally H\"older with H\"older exponent $\gamma$.
Moreover, for each $\delta > 0$, the H\"older seminorm of $u$ in $\Omega^\delta$ depends on the data 
and $\mathrm{osc}_{\Omega_\nu}(u) \delta^{-\gamma}$.
\end{teo}


\subsection{Extension to L\'evy-Ito Operators.} \label{sec:reglevyito}
We present an important extension of our regularity results over equations associated to nonlocal operators 
in \textsl{L\'evy-Ito form}: for $x \in \R^N$ and a bounded function $\phi \in C^2(\bar{B}_\delta(x))$ for some $\delta > 0$, 
we consider $I_x^j$ defined as
\begin{equation}\label{operatorLI}
I_x^j(u, x) = \int_{\R^N} [u(x + j(x,z)) - u(x) - \mathbf{1}_B \langle Du(x), j(x,z) \rangle]\nu(dz), 
\end{equation}
where $\nu$ is a positive regular measure in $\R^N$. The function $j : \R^N \times \R^N \to \R^N$ should be understood as a \textsl{jump function}, 
whose basic assumption concerns the following bound for the jumps, which is  uniform in $x$.

\medskip
\noindent
{\bf (J1)} There exists $C_j > 0$ such that, for all $x \in \R^N$
\begin{equation*}
|j(x, z)| \leq C_j |z|.
\end{equation*}

\medskip

We remark that given $\nu$ and $j$ as above, it is possible to define the associated $x$-dependent measure $\nu_x^j$ 
as the push forward of the measure $\nu$ through the function $j(x, \cdot)$. That is, $\nu_x^j$ is defined as
\begin{equation}\label{pushforward}
\int_{\R^N} f(y) \nu_x^j(dy) = \int_{\R^N} f(j(x,z)) \nu(dz), 
\end{equation}
for each measurable function $f$ satisfying $|f(z)| \leq C\min \{1, |z|^2 \}$ for some $C > 0$. 
It is important to remark that if $\nu$ satisfies (M1),(M2) and $j$ satisfies (J1), then $\{ \nu_x^j \}_{x \in \R^N}$ satisfies 
(M1), (M2) too, where the associated constants now depend on $C_j$.

We also notice that in the case the family of measures $\{ \nu_x^j\}_x$ satisfies (M1)-(M2) with $\sigma \in (0,1)$, then we do not need 
to compensate the integrand and $I_x^j$ is defined as
\begin{equation}\label{operatorLIsigmasmall}
I_x^j(u, x) = \int_{\R^N} [u(x + j(x,z)) - u(x)]\nu(dz).
\end{equation}

For sake of shortness, from this point we mainly argue over $I_x^j$ with the form~\eqref{operatorLI}, but all the results are valid for $I_x^j$
with the form~\eqref{operatorLIsigmasmall} when $\sigma \in (0,1)$.

The following result states the regularity result up to the boundary for L\'evy-Ito problems.
\begin{teo}\label{teoholderIto}
Let $\Omega \subseteq \R^N$ be a bounded domain, $A, b_0 > 0$, $\sigma \in (0,2)$, 
a measure $\nu$ satisfying (M1)-(M2) relative to $\sigma$, and a jump function $j$ satisfying (J1). 
Let $I_x^j$ as in~\eqref{operatorLI},~\eqref{operatorLIsigmasmall} associated to $\nu$ and $j$.
Let $m > \max\{1, \sigma\}$ and $\theta \in [0, m)$.

Then, any bounded viscosity subsolution $u: \R^N \to \R$ to the problem
\begin{equation}\label{eqteoholderIto}
-I_x^j(u,x) + b_0 \ |Du(x)|^m \leq   A d_\Omega(x)^{-\theta}, \quad x \in \Omega 
\end{equation}
is locally H\"older continuous in $\Omega$ with H\"older exponent $\gamma_0$ given in~\eqref{gammaholderbdy}, and H\"older seminorm 
depending on $\Omega$, the data and $\mathrm{osc}_{\Omega_{\nu^j}}(u)$, where $\Omega_{\nu^j}$ is defined as in~\eqref{Omeganu}
relative to the familty of measures $\{ \nu_x^j\}_{x \in \R^N}$ given by~\eqref{pushforward}.

Moreover, if $\Omega$ has a $C^{1,1}$ boundary, then $u$ can be extended as 
a H\"older continuous function to $\bar{\Omega}$ with H\"older exponent $\gamma_0$.
\end{teo}

\noindent
{\bf \textit{Proof:}} This proof follows the lines of Theorem~\ref{teoholderbdy} and therefore we 
provide only a sketch of the proof in order to show how to treat the L\'evy-Ito form. 

\medskip
\noindent
\textsl{1.- Technical lemmas in the L\'evy-Ito context.}
Under the current assumptions, considering $x_0 \in \R^N$, $C_1, C_2, r > 0$ and $\gamma_0$ as in~\eqref{gammaholderbdy}, $w$ defined in~\eqref{defw}
(with $\gamma = \gamma_0$) satisfies the inequality
\begin{equation*}
\sup \limits_{\xi \in B_1(x)} \{ I_\xi^j(w, x) \} \leq C(C_1 + C_2) \varrho^{-\sigma}(x), 
\quad \mbox{for all} \ x \in B_r(x_0) \setminus \{ x_0 \},
\end{equation*}
where $\varrho$ is defined in~\eqref{varrho} and this time the constant $C$ depends also on $C_j$ arising in (J1). 
This is accomplished replacing $\varrho$ by 
\begin{equation*}
\tilde{\varrho}(x) = \min \{ d_0(x), d_r(x) \} /(4C_j), 
\end{equation*}
in the proof of Lemma~\ref{lemaIw}. Once we get this estimate, taking $C_1 > 0$ as in~\eqref{C1holderbdyfinal} (with $C$ now depending on $C_j$)
we conclude
\begin{equation*}
- \sup \limits_{\xi \in B_1(x)} \{ I_\xi^j (w, x) \} + b_0 |Dw(x)|^m \geq A \varrho^{-\theta}(x) 
\quad \mbox{for} \ x \in B_r(x_0) \setminus \{ x_0 \}, 
\end{equation*}
following directly the arguments given in Lemma~\ref{keylemmabdy}.

\medskip
\noindent
\textsl{2.- Censored L\'evy-Ito operators.} Let $u$ be a bounded subsolution to~\eqref{eqteoholderIto}. 
Arguing as in Lemma~\ref{lemacensored}, the L\'evy-Ito analogous to inequality~\eqref{eqauxteoholderbdy} reads as
\begin{equation*}
-I_\Omega^j(u, x) + \frac{b_0}{2}|Du|^m \leq Ad_\Omega(x)^{-\theta} + C(\mathrm{osc}_{\Omega_\nu^j}(u) + 1) d_\Omega^{-\sigma}(x), x \in \Omega,
\end{equation*}
where $C$ depends on $C_j$ and the censored L\'evy-Ito operator $I_\Omega^j$ is defined as
\begin{equation*}
I_\Omega^j(u, x) = \int \limits_{x + j(x,z) \in \Omega} [u(x + j(x,z)) - u(x) - \mathrm{1}_B \langle Du(x), j(x,z) \rangle] \nu(dz).
\end{equation*}

\medskip
\noindent
\textsl{3.- Conclusion.} Once we localize the equation inside $\Omega$, we follow exactly the same lines 
of the proof of Theorem~\ref{teoholderbdy}.
The corresponding inequality~\eqref{C1holderbdyfinal} this time reads as
\begin{equation}\label{C1holderIto}
C_1 \geq C (A^{1/m} + \mathrm{osc}_{\Omega_{\nu}^j}(u) +1), 
\end{equation}
where $C$ depends on $C_j$.
\qed

The immediate consequence of this theorem is the corresponding control of the oscillation. 
Its proof follows the same lines of the one of Corollary~\ref{corosc}
by using the above theorem.
\begin{cor}\label{coroscIto}
Let $\Omega \subset \R^N$ open and bounded, and assume the hypotheses of Theorem~\ref{teoholderIto} hold.
Assume further the nonlocal operator has a censored nature, that is, the family of measures $\{ \nu_x^j\}_{x \in \R^N}$ defined in~\eqref{pushforward}
satisfies the censored condition~\eqref{censoredcondition}. Then, there exists $K > 0$ such that, for each bounded viscosity 
solution of~\eqref{eqteoholderIto} we have
\begin{equation*}
\mathrm{osc}_\Omega(u) \leq K. 
\end{equation*}
\end{cor}

Following the directions given in Theorem~\ref{teoholderIto}, it is possible to provide an interior regularity result
in the flavour of Theorem~\ref{teoholderint}, as well as regularity results for sublinear Hamiltonians in the flavour of Theorems~\ref{teosublinearbdy} 
and~\ref{teosublinearint}, both in the L\'evy-Ito framework. Additionally, we can provide extensions for the mentioned results 
associated to L\'evy-Ito operators when the domain is unbounded (see Remark~\ref{rmkunbounded}). We omit the details.


\section{Well-Posedness for the Cauchy Problem in L\'evy-Ito Form.}
\label{comparisonsection}

The $x$-dependence of the nonlocal term represents a serious difficulty in the statement of 
the comparison principle for integro-differential equations (see~\cite{Barles-Imbert}),
and this comparison principle is a key tool in the study of the large 
time behavior of evolution equations. However, we are able to prove it in 
the interesting case of nonlocal operators in \textsl{L\'evy-Ito} form defined in~\eqref{operatorLI} and~\eqref{operatorLIsigmasmall}.
It is why, from now on,
we consider the Cauchy problem in L\'evy-Ito form
\begin{eqnarray}
\label{pareqIto} \partial_t u(x, t) - I_x^j(u(\cdot, t), x) + H(x, Du(x,t)) &=& 0 \quad \ (x,t) 
\in Q,\\
\label{initialdata} u(\cdot,0) &=& u_0 \quad x \in \R^N,
\end{eqnarray}
where we recall that $Q=\R^N\times (0,\infty).$

We start with the assumptions. Over $\nu$ we require the classical assumption

\medskip
\noindent
{\bf (M)} There exists $C_\nu > 0$ such that
\begin{equation*}
\int \limits_{\R^N} \min \{ 1, |z|^2 \} \nu(dz) \leq C_\nu.
\end{equation*}

We also require the following compatibility condition among $j$ and $\nu$.

\medskip
\noindent
{\bf (J2)} For each $\delta > 0$, there exists $C_\delta > 0$ such that, for each $x, y \in \R^N$ we have
\begin{equation*}
\begin{split}
& \int \limits_{B_\delta} |j(x, z) - j(y, z)|^2 \nu(dz) \leq C_\delta |x - y|^2, \\
& \int \limits_{B \setminus B_\delta} |j(x,z) - j(y,z)| \nu(dz) \leq C_\delta|x - y|. 
\end{split}
\end{equation*}

Concerning the Hamiltonian we assume the following conditions.

\medskip
\noindent
{\bf (H1)} There exists $m > 1$ and moduli of continuity $\zeta_1, \zeta_2$ such that, for all $x, y, p, q \in \R^N$ we have
\begin{eqnarray*}
H(y, p + q) - H(x, p) \leq \zeta_1(|x - y|) (1 + |p|^m) + \zeta_2(|q|) |p|^{m - 1}.
\end{eqnarray*}

\medskip
\noindent
{\bf (H2)} Let $m$ be as in (H1). There exists $A, b_0 > 0$ such that for all $\mu \in (0,1)$ we have
\begin{equation*}
H(x, p) - \mu H(x, \mu^{-1} p) \leq (1 - \mu) \Big{(} b_0 (1 - m) |p|^m + A \Big{)}.
\end{equation*}

\medskip

Note that a measure $\nu$ satisfying (M1)-(M2) satisfies (M). 

Concerning (J1)-(J2), let us give an example. Consider
\begin{equation}\label{jexample}
j(x,z) = g(x) z \quad \mbox{for all} \ x,z \in \R^N.
\end{equation}

If $g: \R^N \to \R$ is bounded then (J1) holds but (J2) may fail.
If, in addition, $g$ is Lipschitz continuous and the measure $|z|\nu(dz)$ 
is finite away the origin, then (J2) holds.

If $m > 1$, assumption (H2) implies~\eqref{H}. 
Examples of Hamiltonians satisfying (H1) and (H2) are provided in 
subsection~\ref{examplesubsection}, see~\eqref{Hexample1},~\eqref{Hexample2}.

\begin{remark}\label{rmksigmasmallLevyIto}
In this section we will argue over nonlocal operators $I_x^j$ with the form~\eqref{operatorLI} (that is, nonlocal operators of
order $\sigma \geq 1$). However, the same arguments can be used to get the results related to $I_x^j$ with the form~\eqref{operatorLIsigmasmall},
replacing (M) by the condition
\begin{equation*}
\int_{\R^N} \min \{ 1, |z|\} \nu(dz) \leq C_\nu < +\infty.
\end{equation*}
\end{remark}

Our comparison principle reads as follows
\begin{prop}\label{parcomparison}
Let $\nu$ be a L\'evy measure satisfying (M), $j$ satisfying (J1) and both satisfying (J2). Let $I_x^j$ defined as in~\eqref{operatorLI}
associated to $\nu$ and $j$. Assume $H$ satisfies (H1),(H2) and $u_0 \in C_b(\R^N)$. 

For each $T > 0$, denote $Q_T = \R^N \times (0,T]$. Let $u, v \in L^\infty(\bar{Q}_T)$ for each $T > 0$ be respective viscosity sub and 
supersolution to~\eqref{pareqIto}-\eqref{initialdata}.
Then, $u \leq v$ in $\bar{Q}$.
\end{prop}

We would like to mention that comparison principles for 
problem~\eqref{pareqIto}-\eqref{initialdata} for the sublinear case (that is $m \leq 1$ in (H1)) are proven in~\cite{Barles-Imbert}
and for this reason we concentrate only in the superlinear case.

The following lemma states the initial condition for viscosity sub and supersolutions is satisfied in the classical sense.
\begin{lema}\label{lemat=0}
Let $I_x^j$ defined in~\eqref{operatorLI} with $\nu$ satisfying (M), $j$ satisfying (J1) and $H$ satisfying (H1). Let 
$u, v$  be respectively a viscosity sub and supersolution to problem~\eqref{pareqIto}-\eqref{initialdata}, satisfying local boundedness in $Q$. Then, 
$u(x,0) \leq u_0(x) \leq v(x,0)$ for all $x \in \R^N$.
\end{lema}

We refer to~\cite{DaLio} for a proof of the corresponding result in the second-order setting. The proof for the current case can be obtained by
adjusting the arguments showed in~\cite{DaLio} to the nonlocal framework. 

We prove Proposition~\ref{parcomparison} in a rather indirect way by using the following lemma, 
which will be also used to prove a version of the Strong Maximum Principle valid for our problem
in Section~\ref{SMPsection}.
\begin{lema}\label{lemau-v}
Let $\sigma \in (0,2)$ and let $I_x^j$ defined in~\eqref{operatorLI} with $\nu$ satisfying (M), $j$ satisfying (J1)
and both satisfying (J2). Assume further that $j(\cdot, z) \in C(\R^N)$ for each $z \in \R^N$. 
Let $H$ satisfying (H1),(H2).
Let $u, v \in L^\infty(\bar{Q}_T)$ for all $T > 0$ be respectively a sub and supersolution to~\eqref{pareq}. 
Then, there exists $\bar{c} > 0$ such that, for each $\mu \in (0,1)$, 
the function 
\begin{equation*}\label{omegalemau-v}
\omega(x,t) := \mu u(x,t) - v(x,t)
\end{equation*}
satisfies, in the viscosity sense, the equation
\begin{equation}\label{eqlemau-v}
\partial_t \omega - I_x^j(\omega(\cdot, t),x) - \bar{c} \frac{\zeta_2(|D\omega|)^m}{(1-\mu)^{m - 1} } 
\leq CA(1 - \mu) \quad \mbox{in}\ Q,
\end{equation}
where $A > 0$ appears in (H2), $\zeta_2$ appears in (H1), $\bar{c} = (m^m b_0^{m-1})^{-1}$ and $C > 0$ is 
an universal constant.
\end{lema}

\noindent
{\bf \textit{Proof:}} We start noting that if $u$ is a viscosity subsolution to~\eqref{pareqIto}, denoting $\bar{u} = \mu u$
we have
\begin{equation}\label{eqbaru}
\partial_t \bar{u} - I_x^j(\bar{u}, x) + \mu H(x, \mu^{-1} D\bar{u}) \leq 0 \quad \mbox{in} \ Q, 
\end{equation}
in the viscosity sense.

Let $(x_0, t_0) \in Q$ and $\phi$ a smooth function such that 
$\omega - \phi$ has a strict maximum point at $(x_0, t_0)$. Let $\epsilon > 0$. Doubling variables we consider the function
\begin{equation*}
\Phi(x,y,s,t) := \bar{u}(x,s) - v(y,t) - \tilde{\phi}(x,y,s,t),
\end{equation*}
where $\tilde{\phi}(x,y,s,t) = \phi(y, t) + \epsilon^{-2}|x - y|^2 + \epsilon^{-2}(s-t)^2$.
By its upper semicontinuity, $\Phi$ attains its maximum over the set
\begin{equation*}
\mathcal{K} := \bar{B}_{2 C_j}(x_0) \times \bar{B}_{2 C_j}(x_0) \times [0,t_0 + 1] \times [0, t_0 + 1] 
\end{equation*}
at a point $(\bar{x}, \bar{y}, \bar{s}, \bar{t})$. Moreover, classical argument in the viscosity theory 
allows us to get that, as $\epsilon \to 0$
\begin{equation}\label{proplemau-v}
\begin{split}
& \bar{x}, \bar{y} \to x_0; \quad \bar{s}, \bar{t} \to t_0; \quad \epsilon^{-2}|\bar{x} - \bar{y}|^2, \
\epsilon^{-2}(\bar{s} - \bar{t})^2 \to 0; \\
& \bar{u}(\bar{x}, \bar{s}) \to \bar{u}(x_0, t_0), \ v(\bar{y}, \bar{t}) \to v(x_0, t_0),
\end{split}
\end{equation}
concluding that for all $\epsilon$ suitably small, 
$\bar{s}, \bar{t}\in (0, t_0+1)$ and $\bar{x}, \bar{y}\in \bar{B}_{2 C_j}(x_0).$
Hence, using that $(x,s) \mapsto \Phi(x, \bar{y}, s, \bar{t})$ 
has a local maximum point at $(\bar{x}, \bar{s})$ and 
$(y, t) \mapsto \Phi(\bar{x}, y, \bar{s}, t)$
has a local minimum point at $(\bar{y}, \bar{t})$, 
we can subtract the viscosity inequality for $v$ at $(\bar{y}, \bar{t})$ to the 
viscosity inequality for $\bar{u}$ (given by~\eqref{eqbaru}) at $(\bar{x}, \bar{s})$ to conclude, for each $\delta' > 0$, the inequality
\begin{equation}\label{testinglemau-v}
\mathcal{A} - I^{\delta'} \leq 0,
\end{equation}
where for $\delta' > 0$ we denote 
\begin{equation*}
\begin{split}
I^{\delta'} = & \ I_{\bar{x}}^j[B_{\delta'}^c](\bar{u}(\cdot, \bar{s}), \bar{x}, \bar{p}) 
- I_{\bar{y}}^j[B_{\delta'}^c](v(\cdot, \bar{t}), \bar{y}, \bar{q}) \\
& \ + I_{\bar{x}}^j[B_{\delta'}](\tilde{\phi}(\cdot, \bar{y}, \bar{s}, \bar{t}), \bar{x}) 
- I_{\bar{y}}^j[B_{\delta'}](-\tilde{\phi}(\bar{x}, \cdot, \bar{s}, \bar{t}), \bar{y}), 
\end{split}
\end{equation*}
with 
\begin{equation*}
\begin{split}
\bar{p} & := D_x \tilde{\phi}(\bar{x}, \bar{y}, \bar{s}, \bar{t}) = 2\epsilon^{-2}(\bar{x} - \bar{y}), \\
\bar{q} & := -D_y \tilde{\phi}(\bar{x}, \bar{y}, \bar{s}, \bar{t}) = \bar{p} - D \phi(\bar{y}, \bar{t}),
\end{split}
\end{equation*}
and
\begin{equation*}
\begin{split}
\mathcal{A} =  (\partial_t \tilde{\phi} - \partial_s \tilde{\phi})(\bar{x}, \bar{y}, \bar{s}, \bar{t}) 
+ \mu H(\bar{x}, \mu^{-1} \bar{p}) - H(\bar{y}, \bar{q}).
\end{split}
\end{equation*}

We estimate each term of the inequality~\eqref{testinglemau-v} to get the result. 
We start with $\mathcal{A}$, noting that taking $\epsilon = \epsilon(\mu)$ small enough, we have 
$$
(1 - \mu) (m -1) b_0 - \zeta_1(|\bar{x} - \bar{y}|) > 0.
$$

Then, from (H1),(H2) we get 
\begin{equation*}
\begin{split}
& \mu H(\bar{x}, \mu^{-1}\bar{p}) - H(\bar{y}, \bar{q}) \\
\geq & \ \mu H(\bar{x}, \mu^{-1} \bar{p}) - H(\bar{x}, \bar{p}) + H(\bar{x}, \bar{p}) - H(\bar{y}, \bar{q}) \\
\geq & \ (1 - \mu) (m - 1) b_0 |\bar{p}|^m  - A (1 - \mu) 
- \zeta_1(|\bar{x} - \bar{y}|)(1 + |\bar{p}|^m) - \zeta_2(|D\phi(\bar{y},\bar{t})|) |\bar{p}|^{m - 1} \\
\geq & \ \inf \limits_{\theta \geq 0} \Big{\{} \Big{(} (1 - \mu) (m - 1) b_0 - \zeta_1(|\bar{x} - \bar{y}|) \Big{)} \theta^{m/(m - 1)} 
- \zeta_2(|D\phi(\bar{y},\bar{t})|) \theta \Big{\}} \\
& \ - A (1 - \mu) - \zeta_1(|\bar{x} - \bar{y}|),
\end{split}
\end{equation*}
that is, denoting $\tilde{c} = (m - 1)^{m - 1}/m^m$, we obtain
\begin{equation*}
\begin{split}
\mu H(\bar{x}, \mu^{-1}\bar{p}) - H(\bar{y}, \bar{q}) \geq & 
- \tilde{c} \frac{\zeta_2(|D\phi(\bar{y}, \bar{t})|)^m}{((1 - \mu) (m - 1) b_0 - \zeta_1(|\bar{x} - \bar{y}|))^{m - 1}} \\
& \ - A (1 - \mu) - \zeta_1(|\bar{x} - \bar{y}|),
\end{split}
\end{equation*}
from which we conclude
\begin{equation}\label{Alemau-v}
\begin{split}
\mathcal{A} \geq & \ \partial_t \phi(\bar{y}, \bar{t})
- \tilde{c} \frac{\zeta_2(|D\phi(\bar{y}, \bar{t})|)^m}{((1 - \mu) (m - 1) b_0 - \zeta_1(|\bar{x} - \bar{y}|))^{m - 1}} \\
& \ - A (1 - \mu) - \zeta_1(|\bar{x} - \bar{y}|).
\end{split}
\end{equation}

Now we address the estimate for $I^{\delta'}$ in~\eqref{testinglemau-v}. Using the smoothness of $\phi$, (M) and (J1) we clearly have
\begin{equation}\label{Ideltalemau-v1}
\begin{split}
& \ I_{\bar{x}}^j[B_{\delta'}](\tilde{\phi}(\cdot, \bar{y}, \bar{s}, \bar{t}), \bar{x}) 
- I_{\bar{y}}^j[B_{\delta'}](-\tilde{\phi}(\bar{x}, \cdot, \bar{s}, \bar{t}), \bar{y}) \\
\leq & \ I^j_{\bar{y}}[B_{\delta'}](\phi(\cdot, \bar{t}), \bar{y}) + \epsilon^{-2} o_{\delta'}(1).
\end{split}
\end{equation}

On the other hand, since $(\bar{x}, \bar{y}, \bar{s}, \bar{t})$ is a maximum point for $\Phi$ 
in $\mathcal{K}$, and since $\bar{x}, \bar{y} \to x_0$ as $\epsilon \to 0$, for all $\epsilon$ small enough, by (J1) we have the inequality
\begin{equation*}
\begin{split}
& \bar{u}(\bar{x} + j(\bar{x}, z), \bar{s}) - v(\bar{y} + j(\bar{y}, z), \bar{t}) - (\bar{u}(\bar{x}, \bar{s}) - v(\bar{y}, \bar{t})) \\
\leq & \ \phi(\bar{y} + j(\bar{y}, z), \bar{t}) - \phi(\bar{y}, \bar{t}) 
+ \epsilon^{-2} (|\bar{x} - \bar{y} + j(\bar{x}, z) - j(\bar{y}, z)|^2 - |\bar{x} - \bar{y}|^2),
\end{split}
\end{equation*}
for each $z \in B_1$. Hence, for each $0 < \delta' < \delta < 1$, using this inequality we conclude that
\begin{equation*}
\begin{split}
& I_{\bar{x}}^j[B_{\delta'}^c](\bar{u}(\cdot, \bar{s}), \bar{x}, \bar{p}) 
- I_{\bar{y}}^j[B_{\delta'}^c](v(\cdot, \bar{t}), \bar{y}, \bar{q})  \\
\leq & \ J^\delta - \int_{B \setminus B_\delta} \langle \bar{p}, j(\bar{x}, z) - j(\bar{y}, z) \rangle \nu(dz)  \\
& + I_{\bar{y}}^j[B_\delta \setminus B_{\delta'}](\phi(\cdot, \bar{t}), \bar{y}) 
+ 2 \epsilon^{-2} \int_{B_\delta \setminus B_{\delta'}} |j(\bar{x}, z) - j(\bar{y}, z)|^2 \nu(dz),
\end{split}
\end{equation*}
where
\begin{equation}\label{Jlemau-v}
\begin{split}
J^\delta = \int_{B_\delta^c} \Big{[} & \bar{u}(\bar{x} + j(\bar{x}, z), \bar{s}) - v(\bar{y} + j(\bar{y}, z), \bar{t}) 
- (\bar{u}(\bar{x}, \bar{s}) - v(\bar{y}, \bar{t})) \\
& \ - \mathbf{1}_B \langle D\phi(\bar{y}, \bar{t}), j(\bar{y},z) \rangle \Big{]} \nu(dz).
\end{split}
\end{equation}

Fixing $\delta > 0$ and using (J2) together with~\eqref{proplemau-v}, we conclude that
\begin{equation*}\label{Ideltalemau-v2}
\begin{split}
& I_{\bar{x}}^j[B_{\delta'}^c](\bar{u}(\cdot, \bar{s}), \bar{x}, \bar{p}) 
- I_{\bar{y}}^j[B_{\delta'}^c](v(\cdot, \bar{t}), \bar{y}, \bar{q})  \\
\leq & \ J^\delta + I_{\bar{y}}^j[B_\delta \setminus B_{\delta'}](\phi(\cdot, \bar{t}), \bar{y}) + C_\delta o_\epsilon(1).
\end{split}
\end{equation*}

Hence, joining the last inequality and~\eqref{Ideltalemau-v1} in the definition of $I^{\delta'}$, we conclude that for all $0 < \delta' < \delta$
\begin{equation*}
I^{\delta'} \leq J^\delta + I_{\bar{y}}^j[B_\delta](\phi(\cdot, \bar{t}), \bar{y}) + C_\delta \ o_\epsilon(1) + \epsilon^{-2} o_{\delta'}(1),
\end{equation*}
with $J^\delta$ defined in~\eqref{Jlemau-v}.
Replacing the last inqueality and~\eqref{Alemau-v} into~\eqref{testinglemau-v}, we conclude that
\begin{equation}\label{ineqlemau-v}
\begin{split}
& \partial_t \phi(\bar{y}, \bar{t}) - I_{\bar{y}}^j[B_\delta](\phi(\cdot, \bar{t}), \bar{y}) - J^\delta 
- \tilde{c} \frac{\zeta_2(|D\phi(\bar{y}, \bar{t})|)^m}{((1 - \mu) (m - 1) b_0 - \zeta_1(|\bar{x} - \bar{y}|))^{m - 1}} \\
\leq & (1 - \mu) A + C_{j, \delta} o_\epsilon(1) + \epsilon^{-2} o_{\delta'}(1) + \zeta_1(|\bar{x} -\bar{y}|). 
\end{split}
\end{equation}

But by (J2), the continuity assumption over $j$, the semicontinuity and boundedness of $\bar{u}, v$ in each $\bar{Q}_T$, by using~\eqref{proplemau-v}
we apply Fatou's Lemma concluding that for each $\delta > 0$ fixed, we get
\begin{equation*}
\limsup \limits_{\epsilon \to 0} J^\delta \leq I_{x_0}^j[B_\delta^c](\omega(\cdot, t_0), x_0, D\phi(x_0, t_0)).
\end{equation*}

Hence, letting $\delta' \to 0$ and $\epsilon \to 0$ in~\eqref{ineqlemau-v}, and recalling~\eqref{proplemau-v} we conclude the desired 
viscosity inequality leading to~\eqref{eqlemau-v}.
\qed

We also require the following
\begin{lema}\label{lemapsibeta}
Let $I_x^j$ defined in~\eqref{operatorLI} with $\nu$ satisfying (M) and $j$ satisfying (J1).
Let $\psi \in C^2_b(\R^d)$ satisfying $||\psi||_{C^2(\R^d)} \leq \Lambda$ for some $\Lambda > 0$. For $\beta > 0$, define 
the function 
\begin{equation}\label{psibeta}
\psi_\beta(x) = \psi(\beta^2 x), \quad x \in \R^N. 
\end{equation}

Then, $\psi_\beta$ satisfies
\begin{equation*}
||D\psi_\beta||_\infty \leq  \Lambda \beta^2, \quad ||D^2 \psi_\beta||_\infty \leq \Lambda \beta^4, 
\quad ||I^j_x(\psi_\beta, \cdot)||_\infty \leq \Lambda o_\beta(1),
\end{equation*}
where $o_\beta(1) \to 0$ as $\beta \to 0$.
\end{lema}

\noindent
{\bf \textit{Proof:}} The estimates concerning $D\psi_\beta, D^2 \psi_\beta$ are direct. Concerning the estimate of the nonlocal term, 
for each $x \in \R^d$ we have
\begin{equation*}
\begin{split}
I^j(\psi_\beta, x) \leq  \ \Lambda \beta^4 \int \limits_{B} |j(x, z)|^2 \nu(dz) + 
\Lambda \beta^2 \int \limits_{B_{1/\beta} \setminus B} |j(x,z)|\nu(dz) + 2 \Lambda \int \limits_{B_{1/\beta}^c} \nu(dz).
\end{split}
\end{equation*}

Hence, using (M) and (J1) in the right-hand side of the last inequality, we get
\begin{equation*}
I_x^j(\psi_\beta, x) \leq C_j^2 C_\nu \Lambda \beta^4  +  C_j \Lambda \beta^2 \int \limits_{B_{1/\beta} \setminus B} |z|\nu(dz) 
+ 2\Lambda \ o_\beta(1).
\end{equation*}

Finally, using that $|z| \leq 1/\beta$ in the integral term of the last 
inequality and applying (M), we conclude the estimate for the nonlocal term.
\qed

Using the last three lemmas we are in position to 
prove the comparison principle for~\eqref{pareqIto}-\eqref{initialdata}.

\medskip

\noindent
{\bf \textit{Proof of Proposition~\ref{parcomparison}:}} Let $T > 0$. We will argue over the finite horizon problem
\begin{equation*}
\left \{ \begin{array}{rll} \partial_t u - I_x^j(u, x) + H(x, Du) & = 0 \quad & \mbox{in} \ Q_T \\ 
u(x, 0) & = u_0(x) \quad & x \in \R^N, \end{array} \right .  
\end{equation*}
from which the general result follows by the fact that $T$ is arbitrary.

We assume by contradiction that
\begin{equation}\label{Mcomparison}
M := \sup \limits_{Q_T} \{ u - v \} > 0.
\end{equation}

Denote $R = 2 (||u||_{L^\infty(\bar{Q}_T)} + ||v||_{L^\infty(\bar{Q_T})})$ 
and consider $\psi \in C^2_b(\R^N)$ a nonnegative function with $\psi = 0$ in $B$, $R \leq \psi \leq 2R$ in $B_2^c$ and 
satisfying  $||D\psi||_\infty, ||D^2\psi||_\infty \leq \Lambda$ for some $\Lambda > 0$. For this function $\psi$ and $\beta > 0$, 
consider $\psi_\beta$ as in~\eqref{psibeta}.

Now, for $\eta, \mu \in (0,1)$, consider the function
\begin{equation*}
\bar{\omega}(x,t) = \mu u(x, t) - v(x,t) - \eta t, \quad (x, t) \in Q.
\end{equation*}

Noting that $\bar{\omega} - \psi_\beta \to u - v$ locally uniform in $\bar{Q}_T$ as $\eta, \beta \to 0$ and $\mu \to 1$, by~\eqref{Mcomparison}
we see that $\bar{\omega} - \psi_\beta$ is strictly positive at some point in $\bar{Q}_T$ for all $\eta, \beta$ close to $0$ and $\mu$ close to $1$. 
Hence, by construction of $\psi_\beta$, $\bar{\omega} - \psi_\beta$ attains its maximum in $\bar{Q}_T$ at some point $(x^*, t^*)$, and by 
Lemma~\ref{lemat=0}, taking $\eta, \beta$ smaller and $\mu$ larger if it is necessary, we have $t^* > 0$ for all such as parameters.
At this point, we fix $\eta > 0$ satisfying the above facts.

Now, by Lemma~\ref{lemau-v}, $\bar{\omega}$ is a viscosity subsolution of
\begin{equation*}
\partial_t \bar{\omega} - I_x^j(\bar{\omega}(\cdot, t),x) - \bar{c} \frac{\zeta_2(|D\bar{\omega}|)^m}{(1-\mu)^{m - 1} } 
\leq CA(1 - \mu) - \eta \quad \mbox{in}\ Q_T,
\end{equation*}
and therefore we can use $\psi_\beta$ as a test function for $\bar{\omega}$ at $(x^*, t^*)$, concluding that
\begin{equation*}
- I_x^j(\psi_\beta,x^*) - \bar{c} \frac{\zeta_2(|D\psi_\beta(x^*)|)^m}{(1-\mu)^{m - 1} } \leq CA(1 - \mu) - \eta.
\end{equation*}

Using Lemma~\ref{lemapsibeta}, we conclude from the above inequality that
\begin{equation*}
-(1 +  \bar{c} (1 - \mu)^{1 - m} )o_\beta(1) \leq CA(1 - \mu) - \eta.
\end{equation*}

Letting $\beta \to 0$ and then $\mu \to 1$, we get the contradiction with the fact that $\eta > 0$.
\qed

As it is classical in the viscosity solution's theory, Proposition~\ref{parcomparison} allows the application of Perron's method to 
conclude the existence. In this task, we introduce the additional asumption

\medskip
\noindent
{\bf (H0)} \textsl{There exists a constant $H_0 > 0$ such that $||H(\cdot, 0)||_\infty \leq H_0$.}

\medskip
This assumption allows us to build sub and supersolutions for~\eqref{pareqIto}. The existence result is the following
\begin{cor}\label{corparcomparison}
Let $I_x^j$ defined as in~\eqref{operatorLI}, with $\nu$ satisfying (M), $j$ satisfying (J1) and both satisfying (J2). 
Assume $H \in C(\R^N \times \R^N)$ satisfies (H0)-(H2).
Let $u_0 \in C_b(\R^N)$. Then, there exists a unique viscosity solution 
$u \in C(\bar{Q}) \cap L^\infty(\bar{Q}_T)$ for all $T > 0$ to problem~\eqref{pareqIto}-\eqref{initialdata}.
\end{cor}

A priori bounds for the solution given in Corollary~\ref{corparcomparison} can be derived from the application of comparison principle.
Using ad-hoc sub and supersolutions, if $u$ is the solution of~\eqref{pareqIto}-\eqref{initialdata}, then
\begin{equation}\label{ubound}
||u(\cdot, t)||_{L^\infty(\R^N)} \leq H_0 \ t + ||u_0||_\infty, 
\end{equation}
which means that for fixed time $t$, the function $x \mapsto u(x, t)$ is globally bounded in $\R^N$.

Similar results can be given for the stationary problem~\eqref{eq} in the L\'evy-Ito setting, namely equations with the form
\begin{equation}\label{eqIto}
\lambda u - I_x^j(u, x) + H(x,Du) = 0 \quad \mbox{in} \ \R^N.
\end{equation}

\begin{prop}\label{comparisonstationary}
Let $\lambda > 0$, $I_x^j$ defined in~\eqref{operatorLI} with $\nu$ satisfying (M), $j$ satisfying (J1) and both satisfying (J2). Assume
$H$ satisfies (H0)-(H2). Let 
$u, v$  be bounded viscosity sub and supersolution to equation~\eqref{eqIto}. Then, $u \leq v$ in $\R^N$. 

Moreover, if in addition we assume (H0), then there exists a unique viscosity solution $u \in C_b(\R^N)$ to equation~\eqref{eqIto}, which satisfies
\begin{equation}\label{uboundstationary}
||u||_\infty \leq \lambda^{-1} H_0.
\end{equation}
\end{prop}


\section{Application to Periodic Equations: Large Time Behavior.}
\label{sectionLTB}

In this section we provide the large time behavior result for the problem~\eqref{pareqIto}-\eqref{initialdata} in the case the data are $\Z^N-$periodic.
Hence, we will argue over the problem
\begin{eqnarray}
\label{pareqtorus} \partial_t u - I_x^j(u(\cdot, t),x) + H(x, Du) &= 0  \quad & \mbox{in} \ \Q := \T^N \times (0,+\infty), \\
\label{initialtorus} u(\cdot,0) &= u_0 \quad & \mbox{in} \ \T^N,
\end{eqnarray}
where $I_x^j$ is a nonlocal operator in L\'evy-Ito form defined in~\eqref{operatorLI}
(replacing $\R^N$ by $\T^N$). Of course, the results obtained in this section can be readily extended to the case
the L\'evy-Ito operator has the form~\eqref{operatorLIsigmasmall}, provided the measure $\nu$ is such that $I_x^j$ 
has order strictly less than 1 (see Remark~\ref{rmksigmasmallLevyIto}).

Since problem~\eqref{pareqtorus}-\eqref{initialtorus} is a particular case of~\eqref{pareqIto}-\eqref{initialdata}, comparison
principle, existence and uniqueness hold for this problem under the conditions on the data given in the statement of Proposition~\ref{parcomparison}.
In particular, for the solution $u$ of~\eqref{pareqtorus}-\eqref{initialtorus} we have the a priori estimate~\eqref{ubound}.


\subsection{Strong Maximum Principle.}
\label{SMPsection}

We need some notation for the statement of the Strong Maximum Principle: let $\nu, j$ in the definition of $I_x^j$ and 
for $x \in \R^N$ we define inductively
\begin{equation*}
X_0(x) = \{ x \}, \quad X_{n + 1}(x) = \bigcup_{\xi \in X_n(x)} \{ \xi + j(\xi, \mathrm{supp} \{ \nu \}) \}, \quad \mbox{for} \ n \in \N,
\end{equation*}
and
\begin{equation}
\mathcal{X}(x) = \overline{\bigcup_{n \in \N} X_n}. 
\end{equation}

The Strong Maximum Principle presented here relies in the nonlocallicity of the operator under the ``iterative covering property''
\begin{equation}\label{supportnu}
\mathcal{X}(x) = \T^N, \quad \mbox{for all} \ x \in \T^N.
\end{equation}

We can provide three interesting examples where this condition clearly holds. Of course,~\eqref{supportnu} depends on both $\nu$ and $j$,
but we mainly focus on the structure of $\nu$ for which this condition is valid, and therefore we assume in the following examples that
$j(x, z) = z$ for all $x, z \in \R^N$. In this context, the most basic example is the case where there exists $r > 0$ such that 
\begin{equation*}\label{ballinsupport}
B_r \subset \mathrm{supp} \{ \nu \}.
\end{equation*}

A second example where the previous property does not hold, but~\eqref{supportnu} remains valid, is
when $\nu$ has the form~\eqref{crossedmeasure}, namely
\begin{equation*}
\nu(dz) = |z_2|^{-(N + \sigma)} dz_2 \otimes \delta_{0}(z_1) dz_1 + |z_1|^{-(N + \sigma)} dz_1 \otimes \delta_{0}(z_2) dz_2, 
\end{equation*}
where $\delta_0$ is the Dirac measure supported at $0$ and $\otimes$ is the measure product. 

The third example strongly takes into account the topology of the torus. In (say) $\T^2$, consider $L \subset \T^2$ a line of irrational slope, 
that is, $L: z_2 = \alpha z_1$, with $\alpha$ irrational. Let $\tilde{\nu}$ be the $1$-dimensional Haussdorff measure in $\T^2$ and 
let $l \subset L$ with $\tilde{\nu}(l) > 0$. Then, the measure $\nu = \mathbf{1}_{l}(z) \tilde{\nu}(dz)$ satisfies the 
assumption~\eqref{supportnu}.

The strong maximum principle is stated through the following
\begin{prop}\label{SMP}
Let $\sigma \in (0,2)$ and let $I_x^j$ defined in~\eqref{operatorLI} with $\nu$ satisfying (M), $j$ satisfying (J1) 
with $j(\cdot, z) \in C(\T^N)$ for each $z \in \R^N$, and $\nu, j$ satisfying (J2) and~\eqref{supportnu}.
Consider $H$ satisfying (H0)-(H2), with $\zeta_2$ in (H1) such that 
$\zeta_2(s) = c|s|$ for some $c > 0$.
Let $u$ be a $\Z^N$-periodic viscosity subsolution to~\eqref{pareqtorus}, and $v$ a $\Z^N$-periodic viscosity supersolution to~\eqref{pareqtorus}, 
such that there exists $(x_0, t_0) \in \Q$ satisfying 
\begin{equation*}
(u - v)(x_0, t_0) = \sup \limits_{\Q} \{ u - v\}.
\end{equation*}

Then, the function $u - v$ is constant in $\mathbb{T}^n \times [0, t_0]$. Moreover, we have
\begin{equation*}
(u - v)(x,t) = \sup \limits_{x \in \T^N} \{ u(x, 0) - v(x, 0) \}, \quad \mbox{for all} \ (x,t) \in \bar{\Q}.
\end{equation*}
\end{prop}

The following lemma is a consequence of the comparison principle, see~\cite{Barles-Souganidis}.
\begin{lema}\label{lemaSMP}
Assume assumptions of Proposition~\ref{parcomparison} hold. Let $u, v$ be locally bounded sub and supersolution to equation~\eqref{pareqtorus}
and for $t \in [0,+\infty)$, define 
\begin{equation}\label{m}
\kappa(t) = \sup \limits_{x \in \T^N} \{ u(x, t) - v(x, t)\}. 
\end{equation}

Then, for all $0 \leq s \leq t$, we have $\kappa(t) \leq \kappa(s)$.
\end{lema}

Now we are in position to prove the strong maximum principle.

\medskip
\noindent
{\bf \textit{Proof of Propostion~\ref{SMP}:}} We divide the proof in several parts.

\medskip
\noindent
\textsl{1.- Preliminaries.}
Under the definition of $\kappa$ in~\eqref{m}, we must prove that for each $(x, t) \in \T^N \times [0,t_0]$
\begin{equation*}
(u - v)(x,t) = \kappa(0).
\end{equation*}

However, since $\kappa(t_0)$ is a global maximum value of $\kappa$ in $[0,+\infty)$, by Lemma~\ref{lemaSMP} we have
$\kappa(t) = \kappa(0)$ for all $t \in [0,t_0]$. Hence, it is sufficient to prove that for each $\tau \in (0,t_0)$ we have
\begin{equation*}
u(x,\tau) - v(x,\tau) = \kappa(\tau), \quad \mbox{for all} \ x \in \T^N,
\end{equation*}
which implies the result up to $\tau = 0$ and $\tau = t_0$ by upper-semicontinuity.

We fix $\tau \in (0, t_0)$ and define the set
\begin{equation*}
\mathcal{M}_\tau = \{ x \in \T^N : (u - v)(x, \tau) = \kappa(\tau) \},
\end{equation*}
which is nonempty by upper-semicontinuity of $u - v$. Hence, with the above facts the proof follows by proving that $\mathcal{M}_\tau = \T^N$.

\medskip
\noindent
\textsl{2.- Localization on time  $\tau$.}
For $\eta > 0$ 
we consider the function
\begin{equation*}
(x,t) \mapsto \tilde{W}(x,t) := u(x, t) - v(x,t) - \eta(t - \tau)^2.
\end{equation*}

Note that for each $(x, t) \in \Q$, we have 
\begin{equation*}
\tilde{W}(x,t) \leq \kappa(t) - \eta(t - \tau)^2 \leq \kappa(\tau) = (u - v)(x_1, \tau) = \tilde{W}(x_1, \tau),
\end{equation*}
for some $x_1 \in \mathcal{M}_\tau$, and therefore the supremum of $\tilde{W}$ in $\Q$ is achieved, and  
each such as maximum point has the form $(x, \tau)$ for some $x \in \mathcal{M}_\tau$. Hence, we clearly have
\begin{equation*}\label{MSMPpartial}
\kappa(\tau) = \sup \limits_{(x,t) \in \mathcal{Q}} \tilde{W}(x,t).
\end{equation*}

\medskip
\noindent
\textsl{3.- Localization around a point in $\mathcal{M}_\tau$.}
From this point we fix $x_\tau \in \mathcal{M}_\tau$
and introduce a function $\psi \in C^2_b(\R)$ with $\psi(0) = 0$, $\psi > 0$ in $\R \setminus \{ 0 \}$
and $\psi(x) = 4 R$ if $|x| \geq 1$, with 
$$
R = ||u||_{L^\infty(\T^N \times [0, t_0 + 1])}.
$$

For $\epsilon > 0$, $x \in \T^N$ define $\psi_\epsilon(x) = \psi(|x - x_\tau|/\epsilon)$. We remark 
that $\psi_\epsilon \in C^2_b(\T^N)$, $\psi_\epsilon(x_\tau) = 0$, $\psi_\epsilon > 0$ in $\R^N \setminus \{ x_\tau \}$ 
and for each $\epsilon > 0$ its first and second derivatives are bounded, depending on $\epsilon$. 

We take $0<\mu < 1$, denote $\bar{u} = \mu u$ and $\omega_\mu = \bar{u} - v$ as in Lemma~\ref{lemau-v}, and
consider the function
\begin{equation*}
(x,t) \mapsto W_\mu(x,t) := \omega_\mu(x,t) - \eta|t - \tau|^2 - (1- \mu) \psi_\epsilon(x). 
\end{equation*}

By upper-semicontinuity of $W_\mu$, there exists $(x_\mu, t_\mu) \in \T^N \times [0, t_0 + 1]$ such that
\begin{equation*}\label{maxpointSMPpartial}
W_\mu(x_\mu, t_\mu) = \sup \limits_{\T^N \times [0, t_0 + 1]} W_\mu,
\end{equation*}
and since $W_\mu \to \tilde{W}$ locally uniform on $\bar{\Q}$ as $\mu \to 1$ we have, up to subsequences,
$(x_\mu, t_\mu) \to (x^*, \tau)$ as $\mu \to 1$, where $x^* = x^*(\epsilon) \in \mathcal{M}_\tau$.

In fact, since $(x_\mu, t_\mu)$ is maximum for $W_\mu$, for all $(x, t) \in \T^N \times [0, t_0 + 1]$ we have 
\begin{equation*}
\begin{split}
W_\mu(x_\mu, t_\mu) = & (u - v)(x_\mu, t_\mu) + (\mu - 1)(u + \psi_\epsilon)(x_\mu, t_\mu) - \eta(t_\mu - \tau)^2 \\
\geq & (u - v)(x,t) + (\mu - 1)(u + \psi_\epsilon)(x, t) - \eta(t - \tau)^2. 
\end{split}
\end{equation*}

In particular, taking the point $(x,t) = (x_\tau, \tau)$ in the right-hand side we obtain
\begin{equation}\label{ineqSMP1}
(u - v)(x_\mu, t_\mu) +  (\mu - 1)(u + \psi_\epsilon)(x_\mu, t_\mu) \geq \kappa(\tau) + (\mu - 1)u(x_\tau, \tau).
\end{equation}

Now, since $t_\mu \in [0, t_0 +1]$ for all $\mu$ close to 1, we have 
\begin{equation*}
(u - v)(x_\mu, t_\mu) \leq \kappa(t_\mu) \leq \kappa(\tau), 
\end{equation*}
and replacing this into~\eqref{ineqSMP1} we get
\begin{equation*}
u(x_\mu, t_\mu) + \psi(|x_\mu - x_\tau|/\epsilon) \leq  u(x_\tau, \tau),
\end{equation*}
that is $\psi(|x_\mu- x_\tau|/\epsilon) \leq 2 R$. By the choice of $\psi$ we conclude that $x_\mu \in B_\epsilon(x_\tau)$ for all $\mu$ close to 1. 
Since $x_\mu \to x^* \in \mathcal{M}_\tau$, we conclude $x^* \in \bar{B}_\epsilon(x_\tau)$.

\medskip
\noindent
\textsl{4.- Using the viscosity inequality for $\omega_\mu$.}
From the above facts, we see that the function $(x,t) \mapsto \phi(x,t) := (1 - \mu)\psi_\epsilon(x) + \eta(t - \tau)^2$ is a test function 
for $\omega_\mu$ at $(x_\mu, t_\mu)$. Then, by Lemma~\ref{lemau-v}, for each $\delta, \epsilon > 0$ we have
\begin{equation*}
\begin{split}
2\eta (t_\mu - \tau) - I_{x_\mu}^j[B_\delta^c] (\omega_\mu(\cdot, t_\mu), x_\mu) - I_{x_\mu}^j[B_\delta]((1 - \mu) \psi_\epsilon, x_\mu) & \\
- \bar{c} (1 - \mu) |D\psi_\epsilon(x_\mu)|^m & \leq CA(1 - \mu),
\end{split}
\end{equation*}
but by (M) and (J1) we have
\begin{equation*}
I_{x_\mu}^j[B_\delta](\psi_\epsilon, x_\mu) \leq C_j |D^2 \psi_\epsilon|_\infty.
\end{equation*}

From this, it follows that
\begin{equation}\label{ineqSMP2}
\begin{split}
2\eta(t_\mu - \tau) - I_{x_\mu}^j[B_\delta^c](\omega_\mu(\cdot, t_\mu), x_\mu) & \\
-(1 - \mu) \Big{(} C_j |D^2\psi_\epsilon|_\infty + \bar{c} |D\psi_\epsilon(x_\mu)|^m + CA \Big{)} & \leq 0.
\end{split}
\end{equation}

Note that for all $\epsilon > 0$, by the smoothness of $\psi_\epsilon$ the term in parenthesis in~\eqref{ineqSMP2} remains bounded as $\mu \to 1$,
meanwhile $t_\mu \to \tau$.
On the other hand, by the continuity of $j$ and (M), by Dominated Convergence Theorem we get
\begin{equation*}
I_{x_\mu}^j[B_\delta^c](\omega_\mu(\cdot, t_\mu), x_\mu) \to I_{x^*}^j[B_\delta^c]((u - v)(\cdot, \tau), x^*) \quad \mbox{as} \ \mu \to 1,
\end{equation*}
where $x^* \in \mathcal{M}_\tau$ is such that $x^* \in \bar{B}_\epsilon(x_\tau)$. Recalling that $(u - v)(x^*, \tau) = \kappa(\tau)$, 
letting $\mu \to 1$ in~\eqref{ineqSMP2} we arrive at
\begin{equation*}
\int_{B_\delta^c} [(u - v)(x^* + j(x^*, z), \tau)  - \kappa(\tau)] \nu(dz) = 0,
\end{equation*}
and since $x^* \in \bar{B}_\epsilon (x_\tau)$, letting $\epsilon \to 0$ we finally conclude
\begin{equation}\label{SMPpasslimit}
\int_{B_\delta^c} [(u - v)(x_\tau + j(x_\tau, z), \tau)  - \kappa(\tau)] \nu(dz) = 0.
\end{equation}

\medskip
\noindent
\textsl{5.- Conclusion.}
Since $\delta > 0$ is arbitrary, we conclude $(u - v)(x,\tau) = \kappa(\tau)$ for all $x \in X_1(x_\tau)$. 
Hence, we can proceed in the same way as above, concluding by induction that $(u - v)(x, \tau) = \kappa(\tau)$ 
for all $x \in \bigcup_{n \in \N} X_n(x_\tau)$. Finally, by upper-semicontinuity of $u - v$ and~\eqref{supportnu} we conclude the result.
\qed

\begin{remark}
In Proposition~\ref{SMP}, the assumption on the continuity of $j$ can be dropped. For instance, it is used to pass to the limit 
in~\eqref{SMPpasslimit}. In this direction, note that if $g \in C(\T^N)$ we can write
\begin{equation*}
\begin{split}
|g(x^* + j(x^*, z)) - g(x_\tau + j(x_\tau, z))| \leq \zeta(x^* + j(x^*, z) - x_\tau - j(x_\tau, z)),
\end{split}
\end{equation*}
where $\zeta$ is the modulus of continuity of $g$. However, it is known  
that a modulus of continuity may be assumed to satisfy 
that $\zeta(t) \leq \zeta(\rho) + \rho^{-1} t$ for each $t, \rho > 0$ (see~\cite{Ishii}). Using this, we conclude
\begin{equation*}
|g(x^* + j(x^*, z)) - g(x_\tau + j(x_\tau, z))| \leq \zeta(\rho) + \rho^{-1}(|x^* - x_1| - |j(x^*, z) - j(x_\tau, z)|)  
\end{equation*}
for all $\rho > 0$.  Hence, using (J1) we can make $x^* \to x_\tau$ and then letting $\rho \to 0$
to get the desired convergence without asking continuity on $j$.

Additionally, instead of assuming $\zeta_2(s)=c|s|,$ it is enough
to ask that
\begin{equation*}
\zeta_2(s) s^{(1 - m)/m} \to 0 \quad \mbox{as} \ s \to 0. 
\end{equation*}
\end{remark}


\subsection{The Ergodic Problem.} 
Roughly speaking, solving the \textsl{ergodic problem} means pass to the limit as $\lambda \to 0$ in the stationary periodic problem 
\begin{equation}\label{eqtorus}
\lambda u -I_x^j(u, x) + H(x,Du) = 0 \quad x \in \T^N,
\end{equation}
whose existence and uniqueness for $\lambda > 0$ holds by Proposition~\ref{comparisonstationary}. 
Hence, the required compactness of the family of solutions $\{ u_\lambda \}$ is typically obtained by regularity results
which are independent of $\lambda$. 
\begin{prop}\label{propergodic}
Let $\sigma \in (0,2)$
and $I_x^j$ defined in~\eqref{operatorLI} with $\nu$ satisfying (M1), (M2) associated to $\sigma$,
$j$ satisfying (J1) with $j(\cdot, z) \in C(\T^N)$ for each $z \in \R^N$, and that $\nu, j$ satisfy (J2) and~\eqref{supportnu}.
Assume $H$ satisfies (H0)-(H2), with $m > \max\{ 1, \sigma \}$ in (H1).
Then, there exists a unique constant $c \in \R$  for which the stationary ergodic problem
\begin{equation}\label{ergodic}
-I_x^j(u, x) + H(x,Du) = - c, \quad \mbox{in} \ \T^N 
\end{equation}
has a solution $w \in C^{(m - \sigma)/m}(\T^N)$. Moreover, $w$ is the unique continuous solution of~\eqref{ergodic}, up to an additive constant.
\end{prop}

\noindent
{\bf \textit{Proof:}} Let $\lambda > 0$ and consider the periodic stationary problem~\eqref{eqtorus}. By Proposition~\ref{comparisonstationary} 
we have the existence and uniqueness of a solution $u_\lambda$ to this problem which,  by~\eqref{uboundstationary}, satisfies the estimate
$||u_\lambda||_\infty \leq \lambda^{-1} H_0$.
Thus, by Theorem~\ref{teoholderIto} we show that
$u_\lambda \in C^{(m - \sigma)/m}(\T^N)$ with H\"older seminorm independent of $\lambda$ or $||u_\lambda||_\infty$.

Now, denote $w_\lambda = u_\lambda - u_\lambda(0)$ which satisfies the equation
\begin{equation}
\lambda u -\I_x^j(u, x) + H(x, Du) = - \lambda u_\lambda(0), \quad \mbox{in} \ \T^N. 
\end{equation}

Using Theorem~\ref{teoholderIto} we see that the family
$\{ w_\lambda \}_{\lambda \in (0,1)}$ is uniformly bounded and that this family is equi-H\"older with exponent $(m - \sigma)/m$.
Hence, by Arzela-Ascoli Theorem, there exists $w \in C^{(m - \sigma)/m}(\T^N)$ such that 
$w_\lambda \to w$ as $\lambda \to 0$, uniformly on $\T^N$. Additionally, we have the existence of a constant $c \in \R$ such that 
$\lambda u_\lambda(0) \to c$ as $\lambda \to 0$.
By standard stability results for viscosity solutions (see~\cite{Barles-Imbert},~\cite{Alvarez-Tourin} and~\cite{usersguide}), 
we have the pair $(w, c)$ found above is a (viscosity) solution to~\eqref{ergodic}.

If $(w_i, c_i)$, $i =1,2$ are two solutions for~\eqref{ergodic}, then we see that $v_i(x,t) = w_i(x,t) + c_i t, i=1,2$ are two solutions to the 
Cauchy problem~\eqref{pareqtorus} with initial data $w_i$. Hence, by comparison principle we conclude that
\begin{equation*}
v_1(x,t) - ||w_1 - w_2||_\infty \leq v_2(x,t), \quad \mbox{for all} \ (x,t) \in \Q,
\end{equation*}
and therefore, we obtain $(c_1 - c_2) t \leq 2||w_1 - w_2||_\infty$. Dividing by $t$ and letting $t \to +\infty$ we obtain that
$c_1 \leq c_2$. Exchanging the roles of $w_1$ and $w_2$, we get 
$c_1=c_2=c$ and therefore $c$ is unique. Moreover, for each $t \in [0,+\infty)$ we have
\begin{equation*}
\sup \limits_{x \in \T^N} \{ v_1(x,t) - v_2(x,t) \} = 
\sup \limits_{\Q} \{ v_1 - v_2 \} = \sup \limits_{\T^N} \{ w_1 - w_2 \} =: m,
\end{equation*}
and therefore, by Proposition~\ref{SMP} we conclude that for each $x \in \T^N$
\begin{equation*}
w_1(x) = w_2(x) + m,
\end{equation*}
concluding the proof.
\qed


\subsection{Large Time Behavior.} The main result of this section is the following
\begin{teo}\label{teoLTB}
Assume assumptions of Proposition~\ref{propergodic} hold.
Let $u$ be the unique solution to problem~\eqref{pareqtorus}-\eqref{initialtorus}. Then, there exists a pair
$(w, c)$ solution to~\eqref{ergodic} such that
\begin{equation*}
u(x,t) - c t  - w(x) \to 0, \quad \mbox{as} \ t \to +\infty,
\end{equation*}
uniformly in $\T^N$.
\end{teo}

\noindent
{\bf \textit{Proof:}} Here we follow closely the arguments given in~\cite{Barles-Souganidis},\cite{Tchamba} in the local framework 
and~\cite{Barles-Chasseigne-Ciomaga-Imbert} in the nonlocal one.

We assume first that $u_0 \in C^2(\T^N)$. In this case, by using comparison principle it is possible to prove 
that $u$ is Lipschitz in $t$ (see~\cite{Tchamba}), with Lipschitz constant
\begin{equation*}
C^* = ||-I_x^j(u_0, \cdot) + H(\cdot, Du_0) ||_{L^\infty(\T^N)} < \infty.
\end{equation*}

Now, by recalling that (H2) implies~\eqref{H}, for each $t \in (0,+\infty)$ the function $x \mapsto u(x,t)$ is a viscosity subsolution to the problem
\begin{equation*}
-I_x^j(u, x) + b_0|Du|^m \leq  C^* + H_0,
\end{equation*}
with $H_0$ given by (H0). Using Theorem~\ref{teoholderIto} we conclude the unique solution to $u$ of problem~\eqref{pareqtorus}-\eqref{initialtorus} 
is in $C^{\gamma_0, 1}(\Q)$, with $\gamma_0$ defined in~\eqref{gammaholderbdy}.

Note that $u$ and the function $(x,t) \mapsto w(x) + c t$ are solutions to~\eqref{pareqtorus}. Hence, by comparison principle 
we have
\begin{equation}\label{LTB2}
||u(\cdot, t) - w - c t||_\infty \leq ||u_0 - w||_\infty,
\end{equation}
meanwhile, if we define 
\begin{equation}\label{defm}
\kappa(t) = \max_{\T^N} \{ u(\cdot, t) - w - c t\},
\end{equation}
by Lemma~\ref{lemaSMP} we see that $\kappa$ is nonincreasing. Since in addition it is bounded
there exists $\bar{\kappa} \in \R$ such that $\kappa(t) \to \bar{\kappa}$ as $t \to +\infty$.

Now, define the function $(x,t) \mapsto v(x,t) := u(x,t) - c t$. Using~\eqref{LTB2} we obtain
\begin{equation*}
||v(\cdot, t)||_\infty \leq ||w||_\infty + ||u_0 - w||_\infty, \quad \mbox{for each} \ t \geq 0,
\end{equation*}
and by the fact that the family $\{ v(\cdot, t) \}_t$ is equi-H\"older (with exponent $\gamma_0$), by Arzela-Ascoli we can extract a subsequence
$\{ v(\cdot, t_k) \}_k$ with $t_k \to \infty$ as $k \to \infty$ such that
\begin{equation*}
v(\cdot, t_k) \to \bar{v}, \quad \mbox{uniformly in} \ \T^N \ \mbox{as} \ k \to +\infty.
\end{equation*}

Define $v_k(x,t) = v(x, t + t_k)$. Recalling that $v_k$ is solution to
\begin{equation*}
\left \{ \begin{array}{rll} \partial_t v_k - I_x^j(v_k(\cdot, t), x) + H(x, Dv_k) &= - c \quad & \mbox{in} \ \Q \\
v_k(x,0) & = v(x, t_k) \quad & x \in \T^N, \end{array}  \right .
\end{equation*}
and using comparison principle we conclude $\{ v_k \}_k$ satisfies the inequality
\begin{equation}
||v_k - v_{k'}||_{L^\infty(\Q)} \leq ||v(\cdot, t_k) - v(\cdot, t_{k'})||_\infty, 
\end{equation}
for all $t \geq 0$ and $k, k' \in \N$. Hence, $\{ v_k \}_k$ is an uniformly bounded Cauchy sequence in $C(\Q)$ and therefore, 
up to a subsequence, we conclude $v_k \to \tilde{v}$ in $C(\Q)$ as $k \to \infty$, where $\tilde{v}$ solves
\begin{equation*}
\left \{ \begin{array}{rll} \partial_t \tilde{v} - I_x^j(\tilde{v}(\cdot, t),x) + H(x, D\tilde{v}) &= - c \quad & \mbox{in} \ \Q \\
\tilde{v}(x,0) & = \bar{v} \quad & x \in \T^N, \end{array}  \right . 
\end{equation*}

Using the definition of $\kappa$ given in~\eqref{defm}, for each $t \geq 0$ we obtain
\begin{equation*}
\kappa(t + t_k) = \max_{\T^N} \{ v_k(\cdot, t) - w\},
\end{equation*}
and since $\{ v_k \}_k$ is uniformly convergent, we can pass to the limit as $k \to \infty$ concluding that
\begin{equation*}
\bar{\kappa} = \max_{\T^N} \{ \tilde{v}(\cdot, t) - w\} \quad \mbox{for each} \ t \in [0,+\infty),
\end{equation*}
and applying Proposition~\ref{SMP}, for each $(x,t) \in \Q$ we have
\begin{equation*}
\tilde{v}(x,t) = w(x) + \bar{\kappa},
\end{equation*}
and therefore we have $\bar{v} = w + \bar{\kappa}$ in $\T^N$. This implies that $v(x, t) \to w + \bar{\kappa}$. But by 
using the definition of $v$ we have 
\begin{equation*}
||u(\cdot, t) - c t - w - \bar{\kappa}||_\infty = ||v(\cdot, t) - v - \bar{\kappa}||_\infty \to 0 
\end{equation*}
as $t \to \infty$. Replacing $w$ by $w + \bar{\kappa}$, we conclude the result in the case the initial data is smooth.

The general result for $u_0 \in C(\T^N)$ follows by an approximation argument using a sequence of smooth initial data $u_0^\epsilon$
satisfying $u_0^\epsilon \to u_0$ uniformly in $\T^N$ as $\epsilon \to 0$. We refer to~\cite{Tchamba} for details.
\qed

\bigskip

\noindent {\bf Acknowledgements:} 
G.B. and O.L. are partially supported by the ANR (Agence Nationale de
la Recherche) through ANR WKBHJ (ANR-12-BS01-0020).
S.K. is supported by Grant-in-Aid for Scientific Research (No. 23340028) of Japan Society for the Promotion of Science.
E.T. was partially supported by CONICYT, Grants Capital Humano Avanzado, Cotutela en el Extranjero and Ayuda Realizaci\'on Tesis Doctoral.

\bigskip



\begin{thebibliography}{00}
%
\bibitem{Alvarez-Tourin}
Alvarez, O and Tourin, A. {\em Viscosity Solutions of Nonlinear Integro-Differential Equations} Annales de L'I.H.P., section C, vol.13 (1996),
no. 3, 293-317.
%

\bibitem{Bardi-DaLio}
Bardi, M. and Da Lio, F. {\em On the Strong Maximum Principle for Fully Nonlinear Degenerate Elliptic Equations} Arch. Math (Basel), 73 (4),
276-285, 1999.

\bibitem{Barles-book}
Barles, G. {\em Solutions de Viscosite des Equations de Hamilton-Jacobi} Collection ``Mathematiques et Applications'' de la SIAM, no 17, 
Springer-Verlag (1994).
%
\bibitem{Barles1}
Barles, G. {\em A Short Proof of the $C^{0,\alpha}-$regularity of Viscosity Subsolutions for Superquadratic Viscous Hamilton-Jacobi Equations and
Applications} Nonlinear Analysis 73 (2010) 31-47.



\bibitem{Barles-Chasseigne-Imbert-Ciomaga-lip}
Barles, G., Chasseigne, E., Ciomaga, A. and Imbert, C. {\em Lipschitz Regularity of Solutions for Mixed Integro-Differential Equations.}
J. Diff. Eq., 252 (2012), 6012-6060.


\bibitem{Barles-Chasseigne-Ciomaga-Imbert}
Barles, G., Chasseigne, E., Ciomaga, A. and Imbert, C. {\em Large Time Behavior of Periodic Viscosity Solutions for Uniformly Elliptic 
Integro-Differential Equations} Preprint. 

\bibitem{bcgj} 
G. Barles, E. Chasseigne, C. Georgelin and E.Jakobsen
{\em On Neumann type problems for nonlocal equations set in a half space.}
To appear in Trans. AMS.


\bibitem{Barles-Chasseigne-Imbert}
Barles, G., Chasseigne, E. and Imbert, C. {\em On the Dirichlet Problem for Second Order Elliptic Integro-Differential Equations}
Indiana U. Math. Journal, 2008.


%
%
\bibitem{Barles-Imbert}
Barles, G. and Imbert, C. {\em Second-order Eliptic Integro-Differential Equations: Viscosity Solutions' Theory Revisited.} 
IHP Anal. Non Lin\'eare, Vol. 25 (2008) no. 3, 567-585.

\bibitem{Barles-Souganidis}
Barles, G. and Souganidis, P.E. {\em Space-time Periodic Solutions and Long-Time Behavior of Solutions of Quasilinear Parabolic Equations.}
SIAM J. Math. Anal., 32 (2001), 1311-1323 (electronic).

\bibitem{BBC}
Bogdan, K., Burdzy, K. and Chen Z.-Q. {\em Censored Stable Processes.} Prob. Theory and Rel. Fields, Vol. 127, Issue 1 (2003), pp 89-152.

\bibitem{Caffarelli-Silvestre}
Caffarelli, L. and Silvestre, L. {\em Regularity Theory For Nonlocal Integro-Differential Equations.} Comm. Pure Appl. Math, Vol. 62 
(2009), no. 5, 597-638.
%
%
%
\bibitem{Capuzzo-Dolcetta-Leoni-Porretta}
Capuzzo-Dolcetta, I., Leoni, F. and Porretta, A. {\em H\"older Estimates for Degenerate Elliptic Equations with Coercive Hamiltonians.} Trans. Amer.
Math. Soc. 362 (9) 4511-4536 (2010).

\bibitem{Cardaliaguet-Cannarsa}
Cardaliaguet, P. and Cannarsa, P. {\em H\"older Estimates in Space-Time for Viscosity Solutions of Hamilton-Jacobi Equations} CPAM 63(5): 
590-629 (2010).

\bibitem{Cardaliaguet-Rainer}
Cardaliaguet, P. and Rainer {\em H\"lder Regularity for Viscosity Solutions of Fully Nonlinear, Local or Nonlocal, 
Hamilton-Jacobi Equations with Superquadratic Growth in the Gradient} SIAM J. Control and Optimization 49(2): 555-573 (2011).

\bibitem{Chasseigne}
Chasseigne, E. {\em The Dirichlet problem for some nonlocal diffusion equations.} Differential Integral Equations 20 (2007), no. 12, 1389???1404. 


\bibitem{Ciomaga}
Ciomaga, A. {\em On the Strong Maximum Principle for Second Order Nonlinear Parabolic Integro-Differential Equations} Advances in Diff. Equations.
17 (2012), 635-671.

\bibitem{Coville}
Coville, J. {\em Maximum principles, sliding techniques and applications to nonlocal equations}
Electronic Journal of Differential Equations, Vol. 2007(2007), No. 68, pp. 1???23.

\bibitem{Coville2}
Coville, J. {\em Remarks on the Strong Maximum Principle for Nonlocal Operators} Electron. J. Differential Equations, pp. No. 66, 10, 2008.

\bibitem{usersguide} Crandall, M.G., Ishii H. and Lions, P.-L. {\em User's Guide to Viscosity Solutions of Second Order Partial
Differential Equations.} Bull. Amer. Math. Soc. (N.S.), Vol. 27 (1992), no. 1, 1-67.

\bibitem{DaLio}
Da Lio, F. {\em Comparison Results for Quasilinear Equations in Annular Domains and Applications.} Comm. Partial Diff. Equations, 27 (1 \& 2) 
283-323 (2002).

\bibitem{Hitch}
Di Neza, E., Palatucci, G. and Valdinoci, E. {\em Hitchhiker's Guide to the Fractional Sobolev Spaces}
Bull. Sci. Math., 136, (2012), no. 5, 521--573.


\bibitem{FOT}
Fukushima M. , Oshima Y.  and Takeda M. 
{\em  Dirichlet forms and symmetric Markov processes.} de Grueter Studies in Mathematics 19 (1994)

\bibitem{GQY}
Guan Q.Y. {\em  The integration by part of the regional fractional  Laplacian.} Commun. Math. Phys. 266, 289--329 (2006)

\bibitem{GQY2}
Guan Q.Y.  and Ma Z.M. {\em  Reflected symmetric $\alpha$ stable process  and regional fractional  Laplacian.}
Probab. Theory Relat. Fields 134, 649--694 (2006)


%
%
\bibitem{Ishii}
Ishii, H. {\em Existence and Uniqueness of Solutions of Hamilton-Jacobi Equations.} Funkcialaj Ekvacioj, Vol. 29 (1986) 167-188.

\bibitem{Ishii-Nakamura}
Ishii, H. and Nakamura, G. {\em A Class of Integral Equations and Approximation of p-Laplace Equations.} Calc. Var. (2010), no. 37, 485--522. 

\bibitem{Ishii-Lions}
Ishii, H. and Lions, P.L. {\em Viscosity Solutions of Fully Nonlinear Second-Order Elliptic Partial Differential Equations } J. Differential Equations,
83(1) 26-78, 1990.

\bibitem{NJ:Book}
Jacob N. {\em  Pseudo differential operators and Markov process. Vol III . Markov process and Applications}
Imperial College Press, Princeton, 2005.

\bibitem{KP}
Kim P. {\em  Weak convergence of censored and reflected stable processes.}
Stochastic Process. Appl. 116 (2006), no. 12, 1792-1814.

\bibitem{Beginners-Koike}
Koike, S. {\em A Beginner's Guide to the Theory of Viscosity Solutions.} Math. Soc. of Japan, 2004.

%
%


\bibitem{Lasry-Lions}
Lasry, J.M. and Lions, P.L. {\em Nonlinear elliptic Equations with Singular Boundary Conditions and Stochastic Control with State Constraints.}
Math. Ann. 283, 583-630 (1989).


%
\bibitem{Sayah1}
Sayah, A. {\em \'Equations d'Hamilton-Jacobi du emier Ordre Avec Termes Int\'egro-Diff\'erentiels. I. Unicit\'e des solutions de viscosit\'e.}
Comm. Partial Differential Equations 16 (1991), 1057-1074.

\bibitem{Sayah2}
Sayah, A. {\em \'Equations d'Hamilton-Jacobi du emier Ordre Avec Termes Int\'egro-Diff\'erentiels. II. Existence de solutions de viscosit\'e.}
Comm. Partial Differential Equations 16 (1991), 1075-1093.
%
\bibitem{Tchamba}
Tchamba, T.T. {\em Large Time Behavior of Solutions of Viscous Hamilton-Jacobi Equations with Superquadratic Hamiltonian} Asymptot. Anal. 66 (2010)
161-186.


\bibitem{Topp}
Topp, E. {\em Existence and Uniqueness for Integro-Differential Equations with Dominating Drift Terms.} Preprint.
\end{thebibliography}
\end{document}